\documentclass[12pt]{iopart}

\usepackage{iopams}
\usepackage[font=footnotesize,labelfont=bf]{caption}
\usepackage{graphicx}
\usepackage{subfig}
\usepackage[breaklinks=true,colorlinks=true,linkcolor=blue,urlcolor=blue,citecolor=blue]{hyperref}

\expandafter\let\csname equation*\endcsname\relax

\expandafter\let\csname endequation*\endcsname\relax

\usepackage{amsmath}
\usepackage{bm}
\usepackage{bbm}
\usepackage{siunitx}
\usepackage{float}
\usepackage{setspace}
\usepackage{bigints}
\usepackage{algorithm}
\usepackage{algpseudocode}
\usepackage{tikz}
\usetikzlibrary{spy}
\usepackage{epstopdf}

\usepackage{color}

\usepackage{booktabs}

\DeclareMathOperator\supp{supp}

\DeclareMathOperator*{\argmin}{arg\,min}
\DeclareMathOperator\TV{TV}

\newcommand{\Rbb}{\mathbb{R}}
\newcommand{\bx}{\boldsymbol{x}}
\newcommand{\bz}{\boldsymbol{z}}

\newtheorem{theorem}{Theorem}

\newtheorem{definition}{Definition}

\newenvironment{Proof}{\paragraph{Proof:}}{\hfill$\square$}

\usepackage{mathtools}

\begin{document}

\title[2D LSFM image reconstruction]{Mathematical Modeling for 2D Light--Sheet
Fluorescence Microscopy image reconstruction}

\author{Evelyn Cueva$^1$, Matias Courdurier$^2$, Axel Osses$^1$, Victor
Casta\~neda$^3$, Benjamin Palacios$^4$, Steffen H\"artel$^3$.}

\address{$^1$ Departamento de Ingenier\'ia Matem\'atica and Centro de Modelamiento
Matem\'atico, UMI CNRS 2807, FCFM, Universidad de Chile, Chile}
\address{$^2$ Facultad de Matem\'aticas, Pontificia Universidad Cat\'olica de
Chile, Chile}
\address{$^3$ Laboratory for Scientific Image Analysis SCIANLab, Facultad de
Medicina, Universidad de Chile, Chile}
\address{$^4$Department of Statistics, University of Chicago, US}
\ead{ecueva@dim.uchile.cl}

\begin{abstract}
We study an inverse problem for Light Sheet Fluorescence Microscopy (LSFM), where
the density of fluorescent molecules needs to be reconstructed. Our first step is to present a
mathematical model to describe the measurements obtained by an optic camera during an LSFM experiment.
Two meaningful stages are considered: excitation and fluorescence. We propose a paraxial model to
describe the excitation process which is directly related with the Fermi pencil--beam equation. For
the fluorescence stage, we use the transport equation to describe the transport of photons towards the
detection camera. For the mathematical inverse problem that we obtain after the modeling, we present
a uniqueness result, recasting the problem as the recovery of the initial condition for the heat equation
in $\Rbb \times(0,\infty)$ from measurements in a space--time curve.
Additionally, we present numerical experiments to recover the density of the fluorescent molecules
by discretizing the proposed model and facing this problem as the solution of a large and sparse
linear system. Some iterative and regularized methods are used to achieve this objective.
The results show that solving the inverse problem achieves better reconstructions than the direct
acquisition method that is currently used.
\end{abstract}

\noindent{\it Keywords\/}: LSFM, Fermi pencil--beam equation, radiative transport equation,
backward uniqueness, heat equation.


\section{Introduction}

Modern microscopy techniques allow researchers to observe phenomena on a sub--cellular,
cellular and supra--cellular level. The observation of cells at different scales gives insights of
key biological questions within modern science fostering more and more systematic approaches to
understand the essence of life~\cite{abbott2003cell}. Contemporary microscopy offers wide spectra of
different techniques with distinct advantages and disadvantages. Particularly, fluorescence
microscopy allows biologists to observe live specimens and dynamic processes within a tissue or
specimen. This technique is based on the addition of fluorescent molecules named fluorophores,
which attach to target proteins or cellular structures on a sub--cellular or cellular level
like DNA, membranes, cytoskeleton, or extra cellular matrix~\cite{lakowicz2013principles}.
Fluorophores are excited by photons, usually a laser beam, and fluorescent emission is captured by a
photonic detector or camera. Fluorescence microscopes vary in the excitation procedure, observation
and volumetric resolution. In the last decades, fluorescence microscopy became the standard tool for
in vivo and in toto (whole sample) imaging, however, photo--toxicity, photo--bleaching, out--of--focus
contribution and acquisition speed limit its application.

Particularly, Light Sheet Fluorescence Microscopy (LSFM) is a technique which uses a thin light
sheet (plane) to excite the fluorophores in the focal plane of the detection objective~\cite{olarte2018light}. This
technique has some advantages compared to the regular confocal fluorescence microscopes. Thanks to
the perpendicular excitation through the thin plane, an optical sectioning occurs. This
excitation reduces the out--of--focus contribution, due to the light sheet only excites
fluorophores present in the observed focal plane. The photo--toxicity and photo--bleaching are also
trimmed down (the energy load is reduced from $10^3$ E to
E~\cite{keller2008quantitative,ritter2011cylindrical}), allowing acquisition of specimen in--vivo
for long periods of time. Moreover, the reduced out--of--focus contribution improves the edges and
contrast of the images. Additionally, its acquisition speed can achieve a few seconds for an entire
3D scan and it can observe big specimens (in the size of millimeters/centimeters)
\cite{santi2011light}. Thus,
LSFM is currently one of the preferred techniques to acquire a wide range of applications,
especially for big specimen and long observation times, obtaining a reasonable image contrast for cell
segmentation and time resolution for cell tracking~\cite{girkin2018light}.
Another related LSFM technique is the so--called lattice light--sheet microscopy
where the laser beam consists in a very narrow Bessel type lattice, intended to
capture much smaller spatial scales of
nanometers~\cite{chen2014lattice,reynaud2014guide}. In this study, we will only
consider LSFM with gaussian type laser beams.

\begin{figure}[ht]
 \centering
 \includegraphics[width=\linewidth]{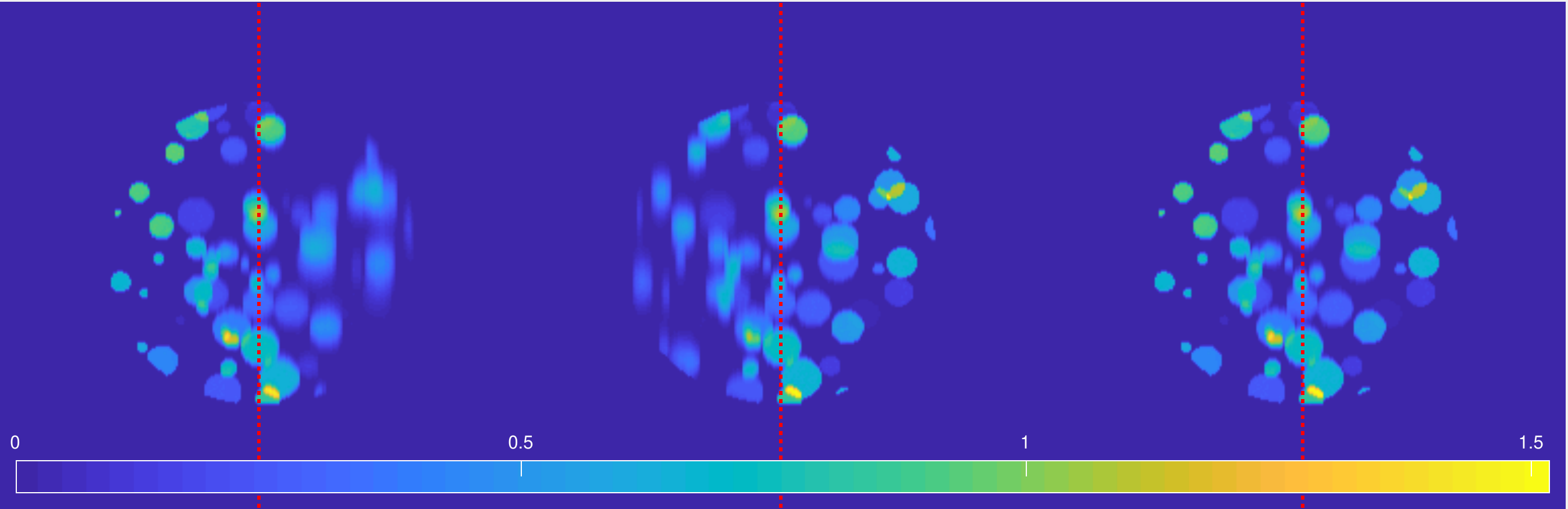}
 \caption{An example of a LSFM image (density of the fluorescent molecules
 $\mu$). The first and second images show the scattering effects observed by
 the camera when left and right excitations are applied. The third one is the
 ``fused image'', taking the best side of the previous ones (as in~\cite{huisken2012slicing, huisken2007even}).}\label{fig:measurements}
\end{figure}

During the image acquisition process, it occurs
that the farther we are from the point of light emission, the higher the loss of
image resolution (see e.g. Figure 2 in~\cite{huisken2012slicing} and Figure 3 in~\cite{huisken2007even}). 
We also see an increasing dominance of blur and shadows as the
laser goes through the object~\cite{huisken2012slicing, reynaud2008light}. The standard reconstruction procedure used to
overcome these issues consists of merging different images by using the opposite
and complementary excitation directions~\cite{huisken2007even, huisken2012slicing, reynaud2008light} (left and right), as in the three
images in Figure~\ref{fig:measurements}. This process is feasible in practice
since the design of the microscope structure is set up in such a way that the
laser beam can illuminate the object from opposite sides preventing the
interference of the lasers.
A critical problem with this merging process is the presence of artifacts in the
middle plane of the final images. On the other hand, there exist calibration
problems in the experimental setting for the acquisition process, such as:
errors in the position and orientation of the lasers respect to the cameras,
object displacements, opposite laser correspondence, etc.

To avoid this merging technique and hence improve the final images, we
establish a mathematical model that allows us to understand the laser behaviour
and the subsequent fluorescence process.  Even more, we propose to study this
imaging technique as an inverse problem, where we seek to reconstruct
the distribution $\mu$ of the fluorophore from the set of (images) measurements
obtained by the camera.

In Section~\ref{sec:model}, we first describe an operator $\mathcal{P}$ that
relates the measurements with our unknown variable $\mu$, identifying two
meaningful stages in an LSFM experiment: excitation and fluorescence. To model
the first stage, we use the Fermi--Eyges pencil--beam equation to describe the
space and angular distributions of the laser beam when it propagates in a
near--transparent object. This equation was first presented by Fermi in 1940 and
studied later by Rossi and Greisen in
\cite[Section~23]{rossi1941cosmic}. In
\cite{borgers1996asymptotic,borgers1996accuracy}, B{\"o}rgers \emph{et al.}
present an asymptotic derivation of the Fermi Pencil--Beam equation from the
Fokker--Planck equation and from the linear Boltzmann equation under two
different conditions.

On the other hand, the fluorescence stage takes place once the
fluorescent molecules has been activated by the laser beam. For the second stage
we use the Radiative Transport Equation (RTE) (see e.g. \cite{bal2009inverse})
to describe how the photons propagate until reaching the collimated
camera.
In this way, we completely define the forward operator $\mathcal{P}$ describing the
proposed mathematical model.

In Section \ref{sec:IP} we summarize the mathematical model obtained and the
description of the inverse problem that we will study.

In Section \ref{sec:uniqueness} we show that there is unique
reconstruction of the function $\mu$ in the proposed inverse problem.
Injectivity of the operator $\mathcal P$ is presented in
Theorem~\ref{th:uniqueness}. We obtain this results by considering the
relationship between the solutions of the Fermi pencil--beam and heat equations.
By interpreting our measurements in terms of heat propagation, we
obtain injectivity of $\mathcal{P}$ by reducing the problem to one of backward
uniqueness for a heat equation from a nontrivial space--time curve, and the
uniqueness for such problem is presented in
Section~\ref{sec:backward_heat}.

Finally, in Sections \ref{sec:numerics} and \ref{sec:numerical_results}, we
present a discretization version of the forward operator to numerically solve
the direct and inverse problems. We propose to find a numerical solution for
the LSFM reconstruction problem by solving a linear system. In this context,
we use different algorithms that are already available to optimally solve this
problems. Mainly, we refer
to~\cite{hansen2010discrete,hansen2018air,gazzola2019ir} where
discrete inverse problems are studied and iterative regulatization methods
for sparse and large--scale problems are detailed.

\section{Mathematical model in LSFM}\label{sec:model}

\subsection{Notation and model scheme}
\begin{figure}[H]
	\centering
	\includegraphics[scale=0.8]{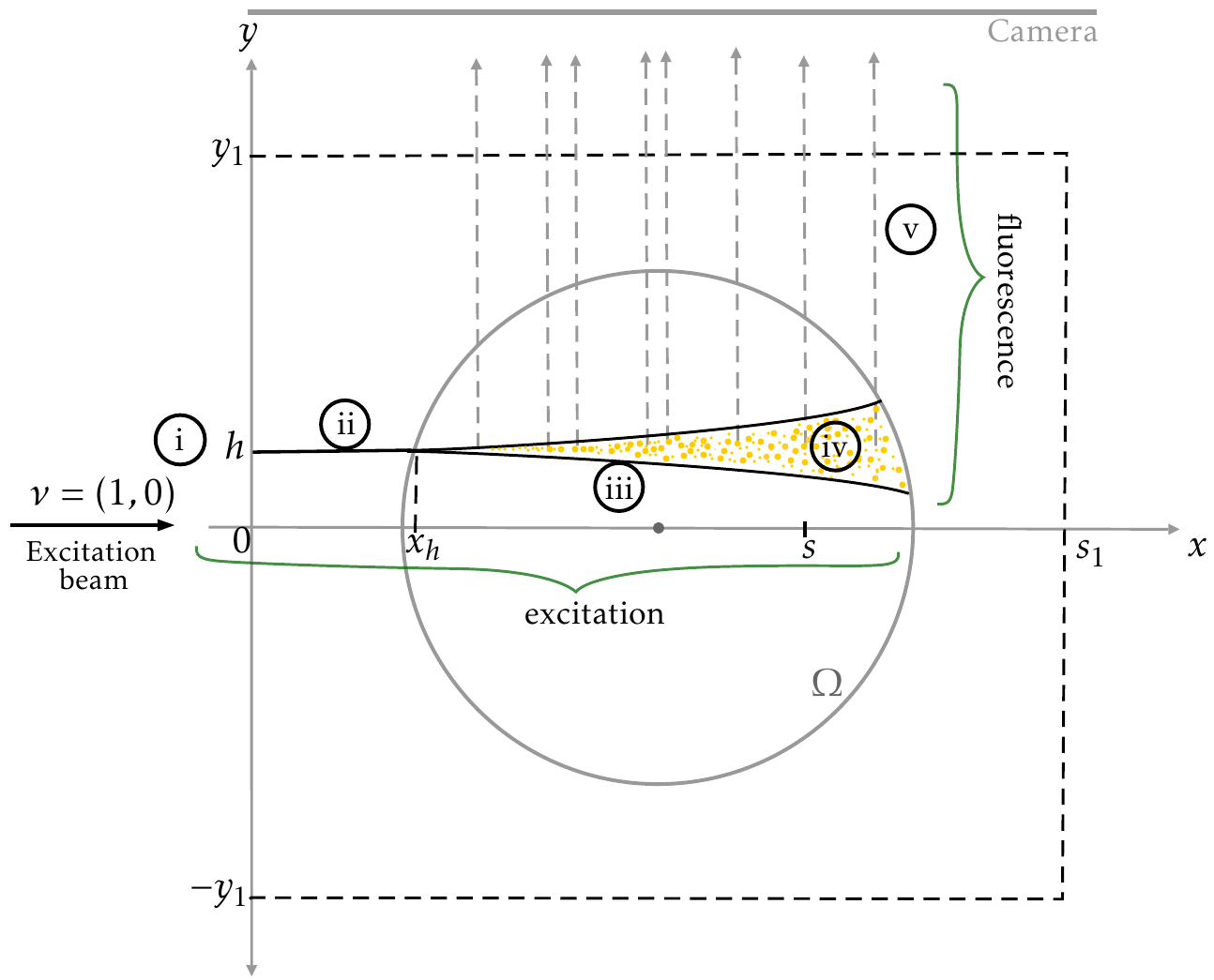}
	\caption{\footnotesize Geometric representation of the excitation and emission beams. An incident laser at point $(0,h)$
    illuminates the object from left and propagates inside the object according to the Fermi pencil--beam 
    equation, exciting the fluorescent molecules within the  sample. Then, the excited fluorescence 
    molecules emits photons in all
	directions. For collimated cameras, only photons emitted in
	straight vertical directions are detected at different positions $s$.\label{fig:Experiment_2D}}
\end{figure}

Let $\Omega\subset\mathbb{R}^2$ be an open set with smooth boundary, which
represents the object studied under the microscope. We assume that
$\Omega$ is contained in the rectangle $[0,s_1]\times [-y_1,y_1]$, for
some $s_1>0,\ y_1 > 0$, both large enough. And for each $h\in[-y_1,y_1]$ we
define $x_h:=\inf\{x:(x,h)\in\Omega\}$ (see in Figure~\ref{fig:Experiment_2D}
the corresponding terms).

The modelling of the LSFM experiment has two main stages, excitation and fluorescence, that are divided in 
the following components (see Figure~\ref{fig:Experiment_2D}):
\begin{enumerate}
\item\label{item:i}  The excitation beam is emitted at the point $ (0,h)$ in the
direction $\nu =(1,0)$. We call $h\in[-y_1, y_1]$ the \emph{height of incidence}.
\item\label{item:ii} The laser follows a free transport equation, without  attenuation or scattering, until 
entering the domain
$\Omega$ at the point $(x_h,h)$.
\item\label{item:iii} Once the laser enters the object, the propagation of the laser is described by the 
Fermi pencil--beam equation (equation \eqref{eq:FermiEq}). We denote by $u := u_h(x, y, \boldsymbol{\omega})$ the intensity of photons at position $(x,y)\in[0,s_1]\times[-y_1,y_1]$ traveling in the direction $ \boldsymbol{\omega}=(\cos(\omega),\sin(\omega))$ for $\omega\in\Rbb/2\pi\mathbb{Z}$. Therefore, the total intensity of excitation photons at $(x,y)$, arising from an incident excitation at $(0,h)$, is  $v_h(x,y)=\int u_h(x,y, \boldsymbol{\omega}) d \boldsymbol{\omega}$.
\item\label{item:iv} The excitation beam reaching $(x,y)$ excites the fluorescent molecules at that point, and the excited fluorophores will be proportional to the density of fluorescent molecules and the excitation intensity. Namely, if $\mu(x,y)$ is the density of fluorescent molecules at $(x,y)$, then the excited fluorophores will be $w_h(x,y) = c~ v_h(x,y)\mu(x,y)$, where $c$ is the activation constant.
\item\label{item:v} The excited fluorescence molecules $w_h$ emit photons in all directions, which propagates according to a linear transport equation (equation \eqref{eq:RTE}). The camera is vertically collimated, hence only measuring those  photons traveling in the direction $(0,1)$. We will denote by $p_h(s)$ the fluorescent measurement at pixel $s\in[0,s_1]$ arising from an excitation at $(0,h)$.
\end{enumerate}

The previous description of LSFM considers some simplifications and does not include all the possible physical phenomena involved in LSFM. The proposed model is a step in trying to understand and tackle difficulties observed in LSFM, like blurring effects among others, and an attempt in trying to improve such imaging technique by analyzing the simplified and related inverse problem.
LSFM can be considered as a particular illumination-detection geometrical setting of Fluorescence Molecular Tomography (FMT) (a review of Fluorescence Molecular Imaging and Fluorescence Molecular Tomography can be found in \cite{ntziachristos2006fluorescence} and \cite{stuker2011fluorescence}), but for a less diffusive media as the one usually considered in FMT. This less diffusive media implies a number of differences between our approach and the usual descriptions used in FMT, namely, in FMT the photon propagation is usually described by a diffusion equation without directionality of photons (see e.g. equation (1) in \cite{lam2005time}, and equations (1) and (2) in~\cite{stuker2011fluorescence}), which translates into a very different mathematical equation for the illumination model. Furthermore, the detection model generally employed in FTM does not allow for directional collimation, and also requires measurements from multiple angles
(see e.g. \cite{ntziachristos2006fluorescence} and \cite{stuker2011fluorescence}).

In the next subsection we present more details about stages \eqref{item:iii} and \eqref{item:v} that we have briefly introduce above.

\subsection{Excitation: the Fermi pencil--beam equation}

In this part we look into the details of stage \eqref{item:iii} above, \emph{i.e.} the propagation of the excitation laser inside the object described by the Fermi pencil--beam equation.

To describe the transport of photons in highly scattering and highly peaked forward regime, a possible model  is the following Fokker--Planck equation (see \cite{bal2009inverse}),
\begin{equation}\label{eq:FokkerPlanckn}
\bomega\cdot \nabla u(\bx, \bomega)+\lambda(\bx,\bomega)u(\bx,\bomega)=\psi(\bx)\Delta_{\bomega}u(\bx,\bomega)
\end{equation}
where, $\bx=(x,y)\in\mathbb R^2$ and $\bomega\in S^1$ is the direction of propagation, with
$\bomega=(\cos(\omega),\sin(\omega))$ for $\omega\in \Rbb/2\pi\mathbb{Z}$.
The quantity $u(\bx,\bomega)$ corresponds to the intensity of photons at the point $\bx$ that are
moving in the direction $\bomega$. The coefficient $\lambda := \lambda_h(\bx, \bomega)$
represents the portion of photons that have been absorbed at the point $\bx$ that were moving in direction $\bomega$. The operator $\Delta_{\bomega} $ is the Laplace--Beltrami operator
on $S^1$ and $\psi(\bx)$ is the diffusion coefficient related to the scattering of the medium. In isotropic media (when $\lambda(\bx, \bomega)=\lambda(\bx)$) and since we are in $\Rbb^2$ (letting $\bomega=(\cos(\omega),\sin(\omega))$), we can rewrite
the Fokker--Planck equation \eqref{eq:FokkerPlanckn} as
\begin{equation}
\label{eq:FokkerPlanck2}
Lu(\bx,\omega)=(\cos(\omega)\partial_x+\sin(\omega)\partial_y +\lambda(\bx)-\psi(\bx)
\partial_\omega^2)u(\bx,\omega)=0.
\end{equation}

And in the case that the diffusion coefficient $\psi(\bx)$ is small enough and the source is spatially  and directionally concentrated, the photons will concentrated along a line and direction determined by the source. Namely, in \cite{borgers1996accuracy} it was shown that under adequate smallness and ellipticity assumptions on the diffusion coefficient, the Fokker--Plank equation
\begin{align*}
L u(x,y,\omega)&=\left(\cos(\omega)\partial_x+\sin(\omega)\partial_y +\lambda(\bx)-
\psi(\bx)\partial_\omega^2\right)u(x,y,\omega)
=0.\\
u(x_h ,y,\omega)&=\delta_h(y)\delta_0(\omega),\qquad x\in (x_h,\infty),y \in\Rbb, \omega\in \Rbb/2\pi\mathbb{Z},
\end{align*}
admits a \emph{paraxial approximation} with $\omega\sim 0$, given by
the Fermi pencil--beam equation
\begin{align}
	\label{eq:FermiEq}
	&{}L_\textnormal{approx} u(x,y,\omega)=\left(\partial_x+\omega\partial_y +\lambda(x,h)-
	\psi(x,h)\partial_\omega^2\right)u(x,y,\omega)=0.\\
	\nonumber &{}u(x_h,y,\omega)=\delta_h(y)\delta_0(\omega)
	,\qquad x\in (x_h,\infty),y \in\Rbb, \omega\in \Rbb,
\end{align}
here we have considered the approximations below inasmuch as $\omega$ is concentrated around zero and satisfies:
\[
\cos(\omega)\approx 1,\quad \sin(\omega)\approx \omega 
\]
and
\[
    |\omega|\ll 1,\ \omega \in \mathbb{R}/2\pi\mathbb{Z}\ \Longleftrightarrow\  |\omega|\ll 1,\   \omega\in \mathbb{R}.
\]
The Fermi equation has been derived from Fokker--Planck in~\cite{borgers1996asymptotic} by means of stereographic--type coordinates on the unit circle and by dropping higher order terms coming from asymptotic expansions with respect to the diffusion magnitude.

Let $\lambda_h(x)=\lambda(x,h)$ and $\psi_h(x)=\psi(x,h)$.
Equation \eqref{eq:FermiEq} can be explicitly solved (see e.g.
\cite{eyges1948multiple}) and the solution for
$x\in (x_h,\infty),\ y \in\Rbb,\ \omega\in \Rbb$ is given by
\begin{equation}\label{eq:paraxial}
u_h(x,y,\omega)
=\exp\left(-\int_{x_h}^x\lambda_h(\tau)d\tau\right)f_Z(\bz),
\end{equation}
where $ \bz =((y-h)-\omega (x-x_h),\omega)^\top$, and where
\[
f_Z(\bz)=\frac{1}{2\pi
\sqrt{\det \Sigma(x,h)}}\cdot\exp\left[-\frac{1}{2}\bz^\top
\Sigma^{-1}(x,h)\bz\right],
\]
with
 \[
 \Sigma(x,h):=\begin{pmatrix}
 E_2&-E_1\\-E_1&E_0
\end{pmatrix}(x,h), \qquad  \Sigma^{-1}(x,h)=\frac 1{\det \Sigma}\begin{pmatrix}
 E_0&E_1\\E_1&E_2
 \end{pmatrix}(x,h),
 \]
and
\begin{align}\label{eq:Ek}
E_k(x,h)=\int_{x_h}^x (\tau-x_h)^k\psi_h(\tau)d\tau,\quad k=0,1,2.
\end{align}
By letting $\Lambda = \begin{pmatrix}
	1 & (x-x_h) \\
	0 & 1 \\
\end{pmatrix} $
(hence $\det(\Lambda)=1$ and
$\Lambda^{-1} = \begin{pmatrix}
	1 & -(x-x_h) \\
	0 & 1 \\
\end{pmatrix} $) then
\[
\bz=\binom{(y-h)-\omega (x-x_h)}{\omega}=
\Lambda^{-1} \binom{y-h}{\omega},
\]
and
\[
f_Z(\bz)=\frac{1}{2\pi
\sqrt{\det \Lambda\Sigma(x,h)\Lambda^\top}}\cdot\exp\left[-\frac{1}{2}\binom{y-h}{\omega}^\top
\left(\left(\Lambda^{-1}\right)^\top\Sigma^{-1}(x,h)\Lambda^{-1}\right)
\binom{y-h}{\omega}\right].
\]
Denoting $\alpha^2=(\Lambda \Sigma \Lambda^\top)_{11} = \left(E_2(x,h)-2(x-x_h)E_1(x,h)+(x-x_h)^2E_0(x,h)\right)$ we get (the marginal distribution on a multivariate normal distribution),
\[
\int_\Rbb f_Z(\bz) dw =\frac{1}{\alpha \sqrt{2\pi}}\exp\left(-\frac{(y-h)^2}{2\alpha^2}\right).
\]
From the solution \eqref{eq:paraxial}, the previous calculation gives us the total excitation intensity at a point $(x,y)\in(x_h,\infty)\times(-y_1,y_1)$ arising from an incident excitation at $(0,h)$, namely
\begin{align}
v_h(x,y)
\nonumber	&=\int_\Rbb u_h(x,y,w) dw = \exp\left(-\int_{x_h}^x\lambda_h(\tau)d\tau\right)\int_\Rbb f_Z(\bz)dw \\
 \label{eq:paraxialxy} &= \frac{1}{\alpha_h(x) \sqrt{2\pi}}\exp\left(-\int_{x_h}^x\lambda_h(\tau)d\tau\right)
	\exp\left(-\frac{(y-h)^2}{2\alpha_h^2(x)}\right),
\end{align}
where
\begin{align}
\nonumber \alpha^2_h(x)&=\left(E_2(x,h)-2(x-x_h)E_1(x,h)+(x-x_h)^2E_0(x,h)\right)\\
\nonumber&=\int_{x_h}^x \psi_h(\tau)[(\tau-x_h)^2-2(x-x_h)(\tau-x_h)+(x-x_h)^2]d\tau\\
\label{eq:alpha_h} &=\int_{x_h}^x (x-\tau)^2\psi_h(\tau)d\tau.
\end{align}

\begin{figure}[ht]
    \begin{center}
        \includegraphics[width=0.8\linewidth]{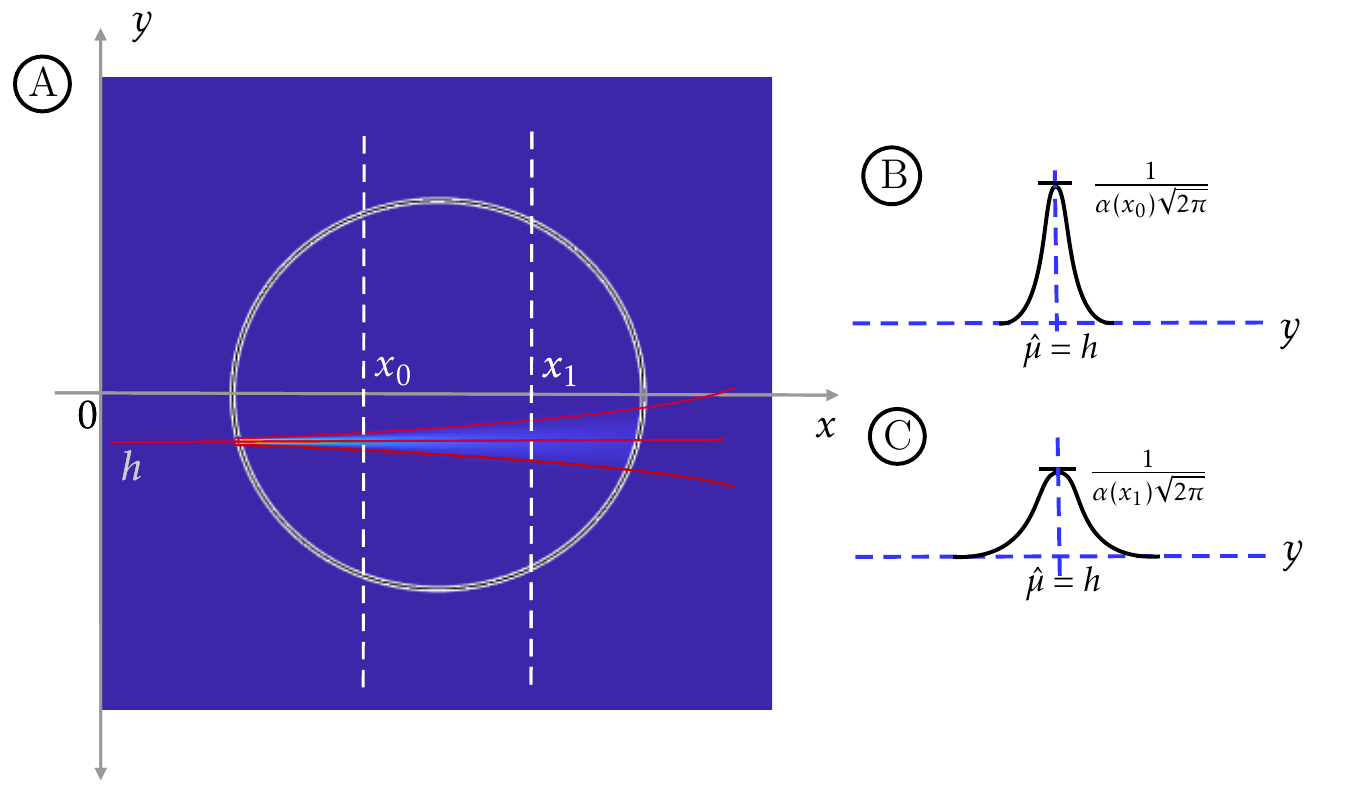}
    \end{center}
    \caption{Graphic interpretation of equation~\eref{eq:paraxialxy}. Figure A shows the function $v_h$ when an illumination is made at height $y=h$. For fix points $x_0$ and $x_1$, the expressions $v_h(x_0,\cdot)$ and $v_h(x_1,\cdot)$ are the density distribution of a normal distribution multiplied by a constant that depends on $\lambda_h$. Figures B and C show these normal distributions. In both cases, the mean is $\hat \mu= h$ with variance $\alpha^2(x_0)$ and $\alpha^2(x_1)$, respectively.}
    \label{fig:vh_gauss}
    \end{figure}

We can notice that for a fix $x$, $v_x(y) = v_h(x, y)$ in~\eref{eq:paraxialxy} is the density function of a univariate normal distribution with mean $h$ and variance $\alpha_h^2(x)$ multiplied by an exponential term depending on $\lambda_h$. This is explained in detail in Figure~\ref{fig:vh_gauss}.

Given the excitation intensity $v_h(x,y)$ and density of fluorescent molecules $\mu(x,y)$, the fluorescent source is $w_h(x,y) = c~v_h(x,y) \mu(x,y)$, and
in the following we provide the details of the model that relates the sources of photons and the measurements obtained at the camera, using the linear transport equation.

\subsection{Fluorescence: Radiative Transfer Equation}

In this detection stage we assume a perfect collimation of the camera in the direction $(0,1)$, this means that only photons travelling parallel to the $y$--axis are measured. The collimation at the camera allows us to remove the positive contribution in the measurements of the scattered photons.

Let us denote by $p_h(\bx,\btheta)$ the intensity of photons at position $\bx\in\Rbb^2$ traveling in a direction $\btheta\in S^1$, arising from an incident excitation at $(0,h)$. We will consider that the propagation of photons is governed by a linear transport equation with attenuation $a$ and source $w_{h}$ (see \cite{bal2011combined, bal2007inverse}), namely we will assume that $p_h$ satisfies,
\begin{align}
\label{eq:RTE}\btheta\cdot \nabla_{\bx}\ p_h(\bx,\btheta )+a(\bx)p_h(\bx,
\btheta\, )=w_{h}(\bx),\quad &{}\forall \bx\in\Rbb^2,\ \btheta\in S^1\\[5pt]
\nonumber\lim_{t\to\infty}p_h(\bx-t\btheta,\btheta\,)=0,\quad &{}\forall \bx\in\Rbb^2,
\ \btheta\in S^1,
\end{align}
where the boundary condition states that there are no external radiation sources, and $w_h$ is supported inside $\Omega$.
Under mild regularity conditions on $w_{h}$ and $a$, the unique solution of equation~\eqref{eq:RTE} is
\[
p_{h}(\bx,\btheta)=\int_{-\infty}^0 w_{h}(\bx+r\btheta)
\exp\left(-\int_{r}^0 a(\bx+\tau\btheta)d\tau\right)dr,
\]
hence providing an expression for the intensity of photons detected at position $\bx$ if collimated in direction $\btheta$.

Since the cameras are outside the bounded object supporting the source, it is useful to consider the total number of photons traveling along lines. In order to do so, let us parametrize the lines in the plane as $L(s,\btheta^\perp)=\{\bx\in \Rbb^2\colon \bx\cdot\btheta =s\}$, where $s\in\Rbb$ is the distance of the line to the origin, $\btheta\in S^1$ is the direction perpendicular to the line, and $\btheta^\perp$, the rotation of $\btheta$ by $\pi/2$, is the direction of the line.
The total intensity of photons along the line $L(s,\btheta^\perp)$ is
 \begin{align}
 p_{h}(s,\btheta^\perp)
 \nonumber&{}=\lim_{\tau\to\infty}p_h(\tau\btheta^\perp+s\btheta,\btheta^\perp)\\
 \nonumber&{}=\int_{\Rbb}w_{h}(r\btheta^\perp+s\btheta)\exp\left(-\int_{r}^\infty
 a(\tau \btheta^\perp +s\btheta)d\tau\right)dr\\
 \label{eq:Measurements} &{}=c\int_{\Rbb}\mu(r\btheta^\perp+
 s\btheta)v_{h}(r\btheta^\perp+s\btheta)\exp\left(-\int_{r}^\infty
 a(\tau \btheta^\perp +s\btheta)d\tau\right)dr,
\end{align}
the last equality is obtained by the assumption $w_h(x,y) = c~ v_h(x,y)\mu(x,y)$ described in \eref{item:iv}.
The Figure~\ref{fig:domain2D} shows an example of the integral along one line.

Under the standard setup of the microscope, the object does not rotate with respect to the camera, hence for the measurements we will consider only the fixed direction $\btheta^\perp=(0,1)$. Rewriting \eqref{eq:Measurements}, and including the expression for $v_h$ given by
\eqref{eq:paraxialxy}, we can finally write an expression for
$p_h(s)=p_h(s,(0,1))$ the intensity of fluorescent photos measured in the camera pixel at position
$s\in[0,s_1]$ arising from an incident excitation at height $h$ (see Figure
\ref{fig:Experiment_2D}):
\begin{align}
\nonumber p_h(s)
&{}=c\int_{\Rbb}\mu(s,r)v_{h}(s,r)\exp\left(-\int_{r}^\infty
a(s,\tau)d\tau\right)dr\\
\label{eq:measure_h}&{}= c\cdot\exp\left(-\int^s_{x_h}\lambda_h(\tau)d\tau\right)\bigintsss_\Rbb
\frac{\mu(s,r)e^{-\int^\infty_r a(s,\tau)d\tau}}{\alpha_h(s)\sqrt{2\pi}}
\exp\left(-\frac{(r-h)^2}{2\alpha_h^2(s)}\right)dr.
\end{align}

We can observe that if $a, \lambda$ and $\psi$ are known, then for
each $h$ fixed, the operator $\mu \mapsto p_{h}(s,\btheta^\perp)$ is a weighted X--ray transform resembling an attenuated X--ray transform with an extra weight. The approach, here presented, considers observations in
multiple heights $h$ for only one angle $\btheta$. But, another interesting problem
can come out if we additionally consider observations for several angles $\btheta\in S^1$, to simultaneously recover $\mu$ and
the attenuation $a$ (or $\lambda$) as in some related works presented in \cite{hertle1988identification,solmon1995identification,stefanov2014identification,courdurier2015simultaneous}.

In the next section, we introduce the measurement operator $\mathcal P$ to study the inverse problem related with the reconstruction of $\mu$ from the expression~\eqref{eq:measure_h}.

\section{Inverse problem}\label{sec:IP}

In this section we will summarize all the elements involved in the description of the measurement operator $\mathcal{P}$, we will discuss about the admissible sections of a domain $\Omega$ where the model $\mathcal P$ is a more adequate description of the phenomena, and we will pose the inverse problems of reconstructing $\mu$ as the inversion of the measurement operator $\mathcal P$.

\subsection{Physical Quantities}

In the previous section we considered the following quantities involved in the phenomena,
\begin{enumerate}
\item $\lambda(x,y)$ describing the attenuation for the incident laser inside the domain.
\item $\psi(x,y)$ describing the diffusion of the laser as it propagates inside the domain.
\item $\mu(x,y)$ the density of fluorescent molecules at each point $(x,y)$ in the domain.
\item $a(x,y)$ describing the attenuation of the fluorescent light inside the domain.
\item $c$ the activation constant, describing the proportion of incident light
that excite the fluorophores.
\end{enumerate}
We will assume $\lambda, \mu, a\in C_\textnormal{pw}(\overline{\Omega})$ and $\psi\in C^1(\overline{\Omega})$, where $C_\textnormal{pw}, C^1$ denote the set of piecewise continuous and continuously differentiable functions, respectively,  we assume that these functions vanish outside of  $\overline{\Omega}$ and that $\psi>0$ in $\overline{\Omega}$. Under these conditions all the solutions to the equations in Section \ref{sec:model} exist and are unique (piecewise continuous regularity could be replaced by $L^1$ regularity).
We recall that we are using the notation $\lambda_h(x):=\lambda(x,h)$
and $\psi_h(x):=\psi(x,h)$.

\subsection{Admissible domain}

It is important to observe that \eqref{eq:paraxial} is a solution to equation \eqref{eq:FermiEq} only under the hypothesis that $\psi_h>0$. Therefore the model for the incident excitation is not as correct after the laser exits the domain $\Omega$, hence equation \eqref{eq:measure_h} describing the fluorescent measurement $p_h(s)$ in pixel $s$ arising from an incident excitation at height $h$, is more adequate if the segment $[x_h,s]\times\{h\}$ is contained in $\overline{\Omega}$. We will consider this aspect for the theoretical part of this work, which motivates the following definitions.

\begin{figure}[H]
 \begin{center}
 \includegraphics[scale=0.8]{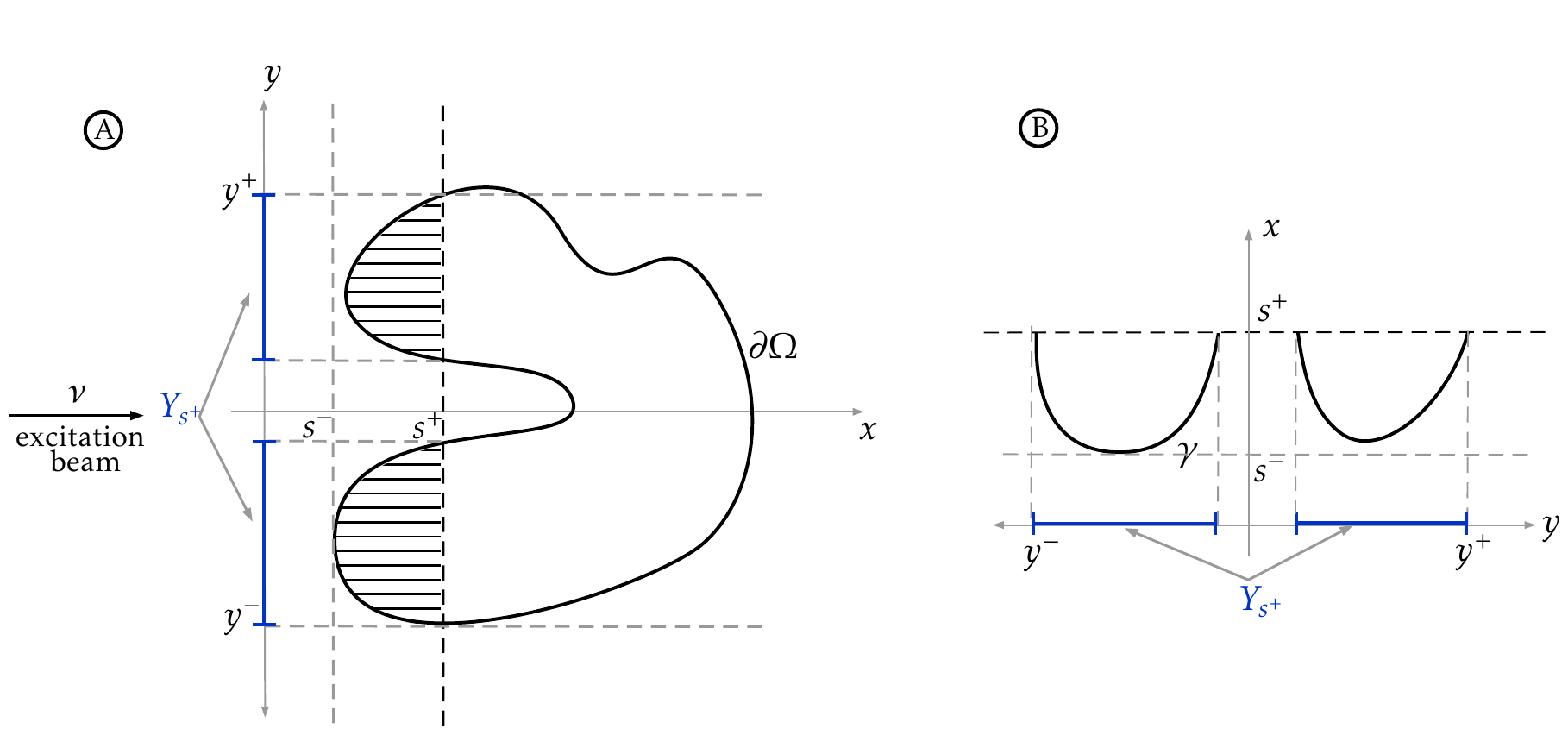}
 \end{center}
 \caption{Example of an admisible domain and the corresponding $\gamma$ function for a generic set $\Omega$. Figure A presents the definition of the quantities $s^-$ and $s^+$
 and the set $ Y_{s^+}$. Figure B shows
 function $\gamma$ and its domain $Y_{s^+}$ in the new coordinates.}
 \label{fig:graph_gamma}
 \end{figure}

\begin{figure}[ht]
\begin{center}
\includegraphics[scale=0.8]{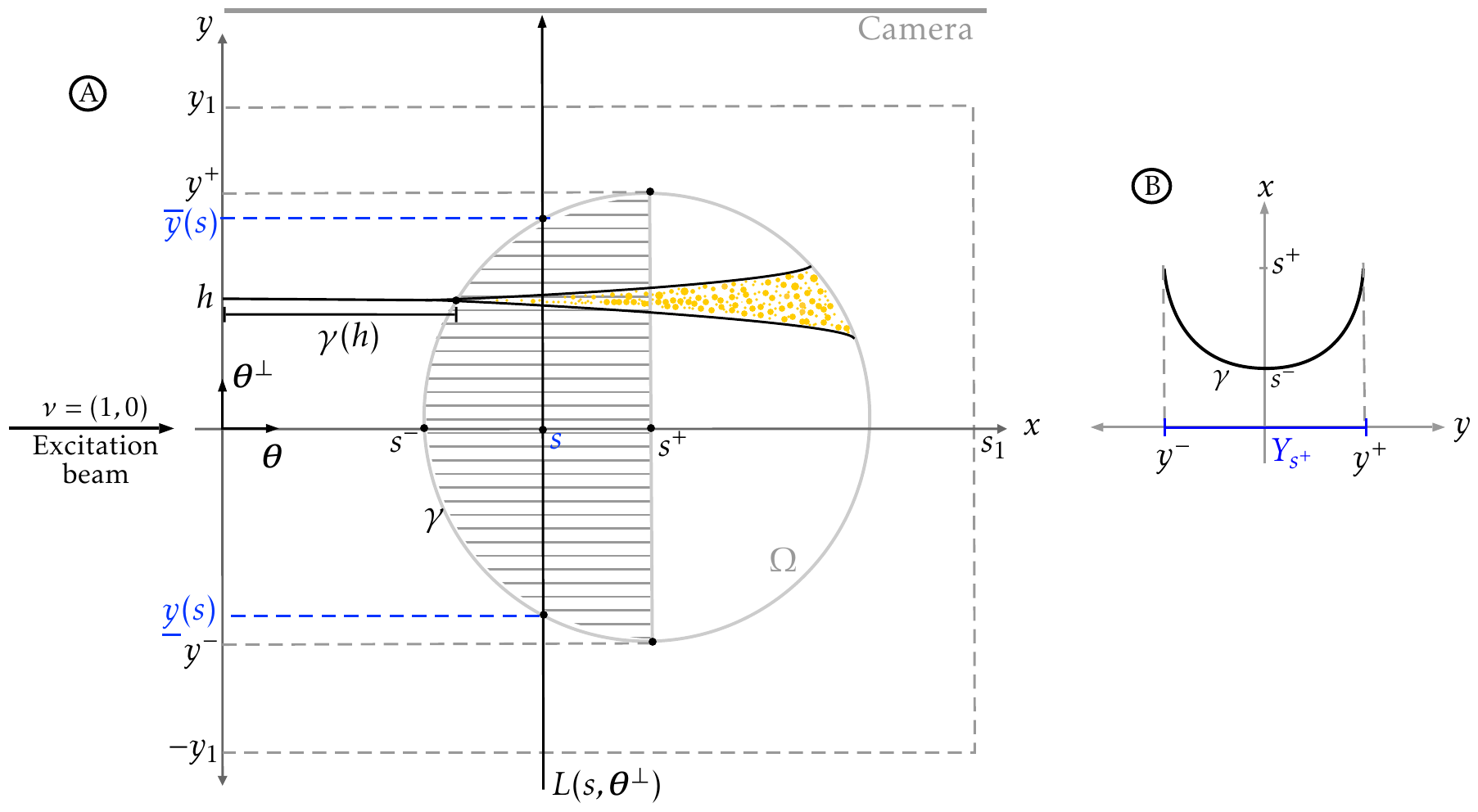}
\end{center}
\caption{Example of an admissible domain for a convex set $\Omega$. On the left side, in figure A, we present its admissible section $\Omega_{\text{ad}}$ filled. All variables
are defined under this scenario. Figure B, at right, shows the corresponding
function $\gamma$ and its domain $Y_{s^+}$.}\label{fig:domain2D}
\end{figure}

\begin{definition}\label{def:admisible_section}
(See Figure \ref{fig:graph_gamma} for an illustration of the following definitions).
Let $\Omega\subset [0,s_1]\times[-y_1,y_1]$ be an open set with smooth boundary. Recall that for $h\in[-y_1,y_1]$ we defined $x_h=\inf\{x:(x,h)\in\Omega\}$. For $s\in[0,s_1]$ define
\begin{align*}
&Y_s=\{h\in[-y_1,y_1] : x_h\leq s\}\\
&s^-=\inf\{s : Y_s\neq \emptyset\},
\end{align*}
and observe that $Y_s\subset Y_r$ for $s < r$.
We say that $s\in [s^-,s_1]$ is admissible if $[x_h,s]\times\{h\}\subset  \overline{\Omega}$, for all $h\in Y_s$. We define $s^+$ as the supremum over the admissible $s$, we define $\underline{y}(s)=\inf(Y_{s})$ and $\overline{y}(s)=\sup(Y_{s})$ for all $s\in[s^-,s^+]$,
and we let $y^-=\underline{y}(s^+), y^+=\overline{y}(s^+)$. We define the admissible section of $\Omega$ as $\Omega_\textnormal{ad}=\{(x,y)\in\Omega : x\leq s^+\}$ and we also define $\gamma:Y_{s^+}\to [0,s^+]$ as $\gamma(h):=x_h$, \emph{i.e.} as the unique smooth function satisfying
\begin{align*}
\Omega_\textnormal{ad} = \{(x,y) : \gamma(y)\leq x\leq s^+\}.
\end{align*}

If the set $\Omega$ is additionally convex, then $Y_{s^+}=[y^-,y^+]$, and if the set $\Omega$ is convex and oriented properly then $\Omega_\textnormal{ad}$ covers half of $\Omega$, in the sense that  at both boundary points $(s^+,y^-)$ and $(s^+,y^+)$ the boundary is tangent to an horizontal line (see Figure \ref{fig:domain2D}).

\end{definition}

Following the discussion above, we will proceed to the theoretical analysis of the inverse problem considering only the admissible section $\Omega_\textnormal{ad}$ of the domain $\Omega$, even though the proposed model could still be used as an approximate description of the whole phenomena in the full domain $\Omega$. Once we are able to solve the inverse problem on an admissible section, the solution to the inverse problem in the full domain follows in a similar fashion as the merging method suggested in~\cite{huisken2007even}. For the right orientation of the camera, which depends on the geometry of the sample, it is possible to solve the inverse problem in $\Omega$ by solving two (or possible more) local problem for admissible regions. This assumes of course the possibility of illuminating the domain from different directions which might be limited by the particular microscope set up. From a numerical point of view, when we leave domain $\Omega$ as we are considering that no diffusion is happening (since $\varphi\in C^1(\Omega)$), integration along lines will be just a rough approximation of the real experiment as in the line shown in Figure~\ref{fig:domain2D}. But if we restrict our analysis to the admissible domain $\Omega_{\textnormal{ad}}$, we guarantee that the integrals along $L(x, \btheta^\perp)$ after excitation at height $y$ with $(x,y)\in \Omega_{\textnormal{ad}}$ fits the exact value given by the model and not just an approximation. We explain this in Figure~\ref{fig:domain2D}.


To complete the framework for the theoretical study we require one more condition with respect to the shape of the domain $\Omega$, prescribed in the following definition.

\begin{definition}
We will say that a domain $\Omega$ is admissible if it satisfies that $\Omega=\Omega_\textnormal{ad}$ and if additionally
$\gamma\in C^1(Y_s)$ and $\gamma'(\underline{y}(s))<0,  \forall s\in(s^-,s^+)$.
\end{definition}

\subsection{Measurements and Inverse Problem}

For the rest of the paper, we will assume that $\Omega=\Omega_\textnormal{ad}$ is an admissible domain,
 in addition to the aforementioned conditions
that $\lambda, \mu, a\in C_\textnormal{pw}(\overline{\Omega}_\textnormal{ad})$, $\psi\in C^1(\overline{\Omega}_\textnormal{ad})$, that these functions vanish at $(x,h)$ if $x<x_h$, and that $\psi>0$ in $\overline{\Omega}_\textnormal{ad}$. In terms of the inverse problem we consider that $\lambda, a$
and $\psi$ are known, while $\mu$ is the unknown quantity.

\begin{definition}[measurement operator]\label{def:measurement_operator}
We define the measurement operator
$\mathcal P$ defined on functions $\mu\in C_\textnormal{pw}(\overline{\Omega}_\textnormal{ad})$
given by (see equation \eqref{eq:measure_h})
\begin{align*}
\mathcal P [\mu] (s,h) = p_h(s), \quad (s, h)\in \Omega_{\text{ad}}.
\end{align*}
\end{definition}
And therefore, the inverse problem consists in recovering $\mu$ from the knowledge of $\mathcal{P}[\mu]$, \emph{i.e.}, we want to study the invertibility of the linear operator $\mathcal P$. 

In next section, we present an injectivity result for the operator $\mathcal P$; this will guarantee that $\ker \mathcal P =\{ 0 \}$ and consequently if the data $p_h(s)$ is in the range of $\mathcal P$, it will uniquely characterize the unknown function $\mu$~\cite{bal2019introduction}. In practice, our measurement operator has to be discretized, and the available data contains noise. Hence, this discretized measurement operator is often not injective, but it will be seen as an approximation of $\mathcal P$, which we will prove is injective. 
 We will overcome the ill-posedness generated by noise data in the discretized inverse problem introducing some regularization techniques as is described in Section~\ref{sec:numerical_results}.

\section{Injectivity of the measurement operator}\label{sec:uniqueness}

For an admissible domain $\Omega_\textnormal{ad}$ and under the hypotheses described in the previous section, we have the following injectivity result for the operator $\mathcal P$.

 \begin{theorem}\label{th:uniqueness}
The measurements $\mathcal{P}[\mu]$ uniquely
 determine the density of fluorophores $\mu$ in $\Omega_{\text{ad}}$., \emph{i.e.} if $\mathcal{P}[\mu](s, h) = \mathcal{P}[\nu](s, h)$ for all
 $(s,h)\in \Omega_{\text{ad}}$ then $\mu(x,y)=\nu(x,y)$ for all $(x,y)\in\Omega_{\text{ad}}$.
 \end{theorem}

This results is a direct consequence of a more localized injectivity property of the linear operator $\mathcal P$, described in the following theorem.

\begin{theorem}
Let $s\in(s^-,s^+)$. If $\mathcal P [\mu] (s,h)=0$ for all $h\in Y_s$ then $\mu(s,y)=0, \forall y\in Y_s$.
\end{theorem}

\begin{Proof}
Let $s\in(s^-,s^+)$ be fixed. Let us recall that for $h\in Y_s$ the measurements take the form (see equations \eqref{eq:alpha_h} and \eqref{eq:measure_h})
 \[
 \mathcal{P}[\mu](s,h) = \exp\left(-\int^s_{\gamma(h)}\lambda_h(\tau)d\tau\right)\bigintsss_\Rbb
 \frac{c\mu(s,r)e^{-\int^\infty_r a(s,\tau)d\tau}}{\sqrt{2\pi \alpha^2_h(s)}}
 \exp\left(-\frac{(r-h)^2}{2\alpha^2_h(s)}\right)dr,
 \]
 where
 \begin{align*}
&\alpha_h^2(s) = \int_{\gamma(h)}^s (s-\tau)^2 \psi(\tau,h)d\tau.
 \end{align*}
 We observe that by letting
\begin{align*}
f(y) &{}:= c\mu(s,y)\exp\left({-\int^\infty_y a(s,\tau)d\tau}\right),\\[8pt]
 g(h) &{}:= \exp\left(\int^s_{\gamma(h)}\lambda_h(\tau)d\tau\right)
\mathcal P[\mu](s,h),\ \text{and}\\[8pt]
\sigma(h) &{}:= \alpha_h^2(s)/2,
\end{align*}
then the theorem reduces to show that  $f(y)=0, \forall y\in Y_s$ whenever $g(h)=0, \forall h\in Y_s$, where
 \begin{equation}\label{eq:data_yb}
 g(h) = \int_\Rbb
 \frac{f(r)}{\sqrt{4\pi \sigma(h)}}
 \exp\left(-\frac{(r-h)^2}{4\sigma(h)}\right)dr.
 \end{equation}
If $U(t,y)$ is the unique solution to the following initial value problem for the heat equation,
\begin{equation}
 \left\{\begin{aligned}
& (\partial_t - \partial^2_y)U(t,y) = 0, \quad(t,y)\in(0,+\infty)\times\mathbb{R},\\
 &U(0,y) = f(y), \textnormal{ if } y \in Y_s,\\
& U(0,y) = 0, \textnormal{ if } y\notin Y_s,  \\
& \lim_{|y|\to\infty} U(t,y) = 0, \forall t>0,
 \end{aligned}
 \right.
 \end{equation}
then
 \[
 U(t,y) = \int_\Rbb \frac{f(r)}{\sqrt{4\pi t}}\exp\left(-\frac{(r-y)^2}{4t}\right)dr,
 \]
 and
\begin{align*}
g(y) = U(\sigma(y),y),\ \forall y\in Y_s,\quad \textnormal{ while } \quad f(y) = U(0,y),\ \forall y\in Y_s.
\end{align*}
Let $\Gamma := \{(\sigma(y),y) : y\in Y_s\}\cup\{(0,y) : y\notin Y_s\}$. Since
$g(y)=0, \forall y\in Y_s$ if and only if $U|_\Gamma =0$, then we can recast our problem as the problem of proving that
\begin{align*}
U|_{\Gamma} = 0 \textnormal{ implies } U(0,y) = 0,\ \forall y \in Y_s.
\end{align*}
This is exactly what Theorem \ref{thm:backward_uniqueness} in the following section shows. But to use Theorem \ref{thm:backward_uniqueness} we need to check that $\Gamma$ satisfies the required conditions, which reduces to prove the following
\begin{enumerate}
\item\label{sigma_i} $\sigma:Y_s\to \Rbb$ is $C^1$.
\item\label{sigma_ii} $\sigma(y)=0$ if $y\in \partial Y_s$.
\item\label{sigma_iii} $\sigma'(y)=0$ whenever $\sigma(y)=0$.
\item\label{sigma_iv} There exists $\delta>0$ such that $\sigma'(y)>0$ for $y\in(\underline{y}(s),\underline{y}(s)+\delta)$.
\end{enumerate}
Let us prove this four points. Recall that for $y\in Y_s$
 \begin{align}\label{eq:sigma}
\sigma(y) = \frac 12\int_{\gamma(y)}^s (s-\tau)^2 \psi(\tau,y)d\tau,
 \end{align}
therefore
 \begin{align}\label{eq:sigma'}
\sigma'(y) = -\frac 12\gamma'(y) (s-\gamma(y))^2 \psi(\gamma(y),y)+
\frac 12\int_{\gamma(y)}^s (s-\tau)^2 \frac{\partial \psi}{\partial y} (\tau,y)d\tau.
 \end{align}
The hypotheses on the regularity of $\gamma$ and $\psi$ clearly imply that
$\sigma\in C^1(Y_s)$ and therefore \eqref{sigma_i} is satisfied. Property \eqref{sigma_ii} follows from the equation \eqref{eq:sigma} and the fact that
if $y\in \partial Y_s$ then $\gamma(y)=s$.
In order to check \eqref{sigma_iii} let us recall that $\psi>0$ in $\Omega$, therefore $\sigma(y)=0$ only if $\gamma(y)=s$ (see equation \eqref{eq:sigma}), in which case equation \eqref{eq:sigma'} implies $\sigma'(y)=0$. To establish \eqref{sigma_iv}, we observe that
if $m=\inf_{(x,y)\in\overline{\Omega}}|\psi (x,y)|>0$ and
$M=\sup_{(x,y)}|\partial \psi/\partial y (x,y)|$ then from equation
\eqref{eq:sigma'}
\begin{align*}
\frac{2\sigma'(y)}{(s-\gamma(y))^2} \geq \Big[-\gamma'(y)m - \frac{1}{3} (s-\gamma(y))M\Big]\stackrel{y\to\underline{y}(s)}{\longrightarrow} -\gamma'(\underline{y}(s))m,
\end{align*}
since $\Omega$ is admissible, $\gamma'(\underline{y}(s))<0$ and
therefore $\sigma'(y)>0$ for $y\in(\underline{y}(s),\underline{y}(s)+\delta]$,
for some $\delta>0$.
 \end{Proof}

 \section{A uniqueness result for the heat equation}\label{sec:backward_heat}

 The purpose of this section is to prove the next result.
 \begin{theorem}\label{thm:backward_uniqueness}
 Let $\sigma(y)\in C^1_c(\mathbb{R})$ and denote $\Gamma =\{(t,y)\in\mathbb{R}^2: t=\sigma(y)\}$. Let
 $$\underline{y}=\inf(\supp\sigma),\quad\overline{y}=\sup(\supp\sigma)$$
 and assume there is $\delta>0$ so that $\sigma'(y)> 0$ in $(\underline{y},\underline{y}+\delta)$.
 If $U(t,y)$ is a solution to the heat equation
\begin{align*}
(\partial_t - \partial^2_y)U(t,y) = 0,\quad (t,y)\in(0,+\infty)\times\mathbb{R},\\
U(t,y)\to 0\;\text{as}\; |y|\to\infty,\ \forall t>0,
\end{align*}
satisfying $\supp U|_{t=0}\subset\supp\sigma$ and $
 U|_{\Gamma} = 0$,  then $U= 0$ everywhere in $(0,+\infty)\times\mathbb{R}$. In particular
$U(0,y)=\lim_{t\to 0^+}U(t,y)=0, \forall y\in\Rbb$.
 %
 \end{theorem}
 \begin{Proof}
Let $T=\sigma(\underline{y}+\delta)$, by hypothesis the restriction of $\sigma$ to the interval $(\underline{y},\underline{y}+\delta)$ has
 an inverse $\rho(t) = \sigma^{-1}(t)\in C^1(0,T)\cap C[0,T]$, and since $\sigma(\underline{y}) = 0$ then $\rho(0) = \underline{y}$. Then we can parameterize the section of $\Gamma$ immediately to the right of $(\underline{y},0)$ as $\{(\rho(t),t):0\leq t\leq T\}$ (see Figure~\ref{fig:Gamma}).
 \begin{figure}
 	\centering
 	\includegraphics[scale=0.9]{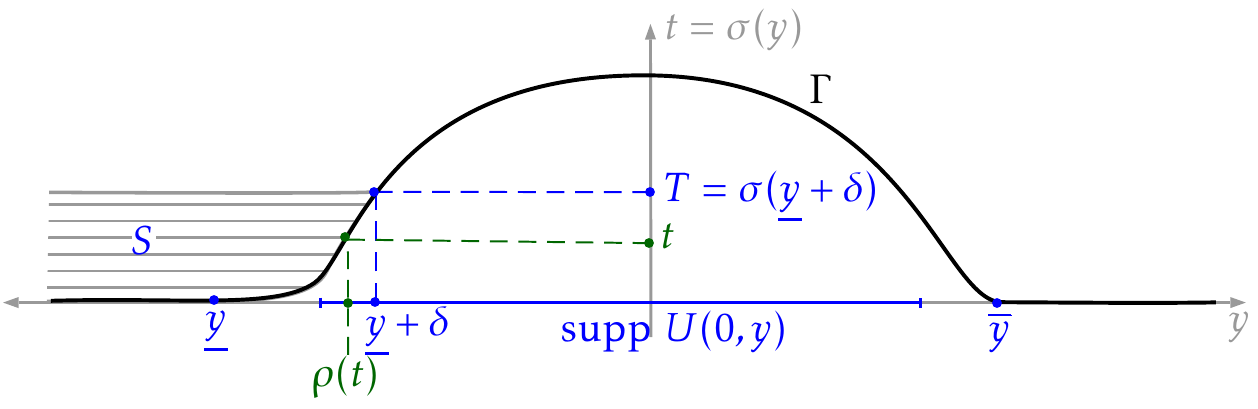}
 	\caption{Curve $\Gamma$ in variables $(y,t)\in \Rbb \times (0,+\infty)$. The
	filled zone $\{(y,t)\colon 0\leq t< T,\ y<\rho(t)\}$ has been denoted by $S$. The assumption that $\supp U|_{t=0}\subset\supp\sigma$ is also represented, since
	$\underline{y}=\inf(\supp\sigma)$ and $\overline{y}=\sup(\supp\sigma)$.}
 	\label{fig:Gamma}
 \end{figure}
 Let us define the following one--sided exterior energy
 $$I(t):= \frac{1}{2}\int^{\rho(t)}_{-\infty}|U(t,y)|^2dy,\quad t\in[0,T),$$
 and notice that for all $t\in(0,T)$
 $$\frac{d}{dt}I(t) = \frac{1}{2}|U(\rho(t),t)|^2\frac{d}{dt}\rho(t) + \int^{\rho(t)}_{-\infty}U(t,y)\partial_tU(t,y)dy ,$$
and the first term in the sum vanishes since $U|_{\Gamma}=0$. On the other hand, since $U$ solves the heat equation and integrating by parts,
 $$
 \begin{aligned}
 \int^{\rho(t)}_{-\infty}U(t,y)\partial_tU(t,y)dy &=  \int^{\rho(t)}_{-\infty}U(t,y)\partial_y^2
 U(t,y)dy\\
 &= U(t,\cdot)\partial_y U(t,\cdot)\Big|^{\rho(t)}_{-\infty} - \int^{\rho(t)}_{-\infty}|\partial_y
  U(t,y)|^2dy,
 \end{aligned}
 $$
and again the first term in the sum vanishes since $U|_{\Gamma}=0$.
Therefore
 $$\frac{d}{dt}I(t) = - \int^{\rho(t)}_{-\infty}|\partial_y U(t,y)|^2dy\leq 0,\quad
 \forall t\in[0,T),$$
and $I(t)$ is a nonnegative decreasing function.
But $\supp U(0,y)\subset\supp\sigma$, implying that $I(0) = 0$ and concluding that $I(t)=0$ for all  $t\in [0,T)$. It follows that
 $$U(t,y)= 0,\quad  \forall t\in[0,T), \forall y<\rho(t),$$
and from classical unique continuation results for parabolic equations (see for instance \cite{lin1990uniqueness}) we deduce that $U$ must vanish in the whole upper-half plane.

 \end{Proof}

 In the next sections, we present the numerical implementation of the direct and inverse problems.

 \section{Discrete direct and inverse problems}\label{sec:numerics}

 The main objective of this and next sections is to present a numerical analysis and solution of the direct and inverse problems. This will allow us to bear out that the diffusion and artifacts, observed during the traditional acquisition process, can be described by the proposed model.

 \subsection{Direct model}

 Here, we present how to simulate our data set using the proposed forward operator $\mathcal P$. Given the fluorescence density $\mu$ in a given domain $\Omega$, we are able to compute the value of $p_h(s)$ for all $s$ thanks to the expression~\eqref{eq:measure_h}.

\begin{figure}[ht]
  \centering
  \includegraphics[width=\linewidth,height=0.35\linewidth]{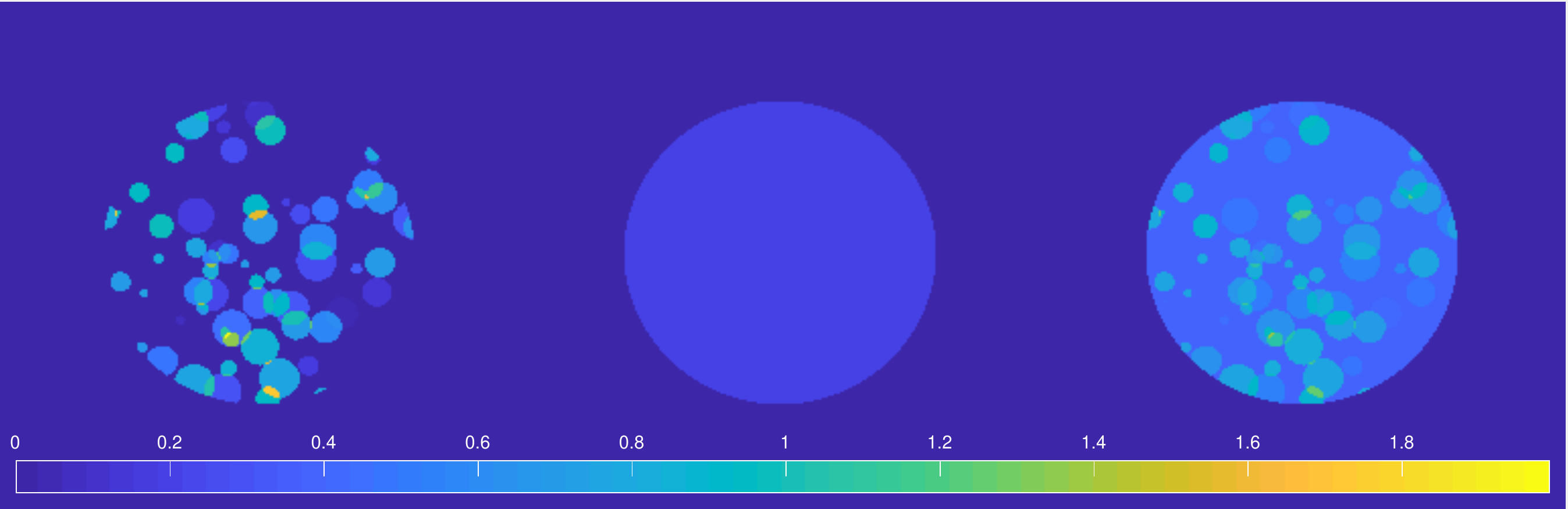}
  \caption{\textbf{From left to right:} fluorophore density distribution ($\mu$), constant and variable attenuation ($\lambda$) for the excitation stage.}
  \label{fig:constant_data}
 \end{figure}

 The density of fluophores $\mu$ and the two cases of attenuation $\lambda$ that we will consider in the experiments are presented in Figure~\ref{fig:constant_data}. The variable attenuation is proportional to the fluorophore density plus a constant value which represents the medium where the object is submerged. We assume that the attenuation of the fluorescence stage $a$ satisfies the relation $a = \hat c \cdot \lambda$. We choose a parameter $\hat c$ so that the diffusion effect got in the numerical experiments remains close to the one observed in the real data. Here, we also assume that the diffusion term $\psi_h$ is proportional to the attenuation $\lambda_h$,
 \emph{i.e.} $\psi_h = \tilde c\cdot \lambda_h$. For all the experiments we set this constant in $\tilde c = 0.6$. Additionally, recalling that $w_h = c\cdot \mu \cdot v_h$, represents the amount of fluorescent molecules that is activated after the excitation process, we took $c=1$ throughout the experiments.

 For all experiments, we work over the domain $\Omega=[0,2]\times[-1,1]$ and with images of size $N \times N$ with $N=257$. The discretization step is given by $\tau = 2/(N+1)$ in $x$ and $y$ axes.
 We start by calculating the values $v_h(x,y)$ over $\Omega$ for a discretized set of excitations points along the interval $[-1,1]$. We take $N$ heights of excitations with step size $\tau$. The excitation points are considered in two directions: left and right, since the support of our object is a circle (as shown in Figure~\ref{fig:constant_data}) by the Definition~\ref{def:admisible_section}, two directions are needed to guarantee the uniqueness of our solution in the whole domain. Then the total amount of excitation points is $2N$.

 The discretization of equation~\eqref{eq:paraxialxy} is straightforward if we approximate the integrals of $\lambda_h$ as finite sums of its pixel intensities, since we are representing $\lambda_h$ as an image of size $N\times N$. The same is considered for the integrals of $\psi_h$ in expression~\eqref{eq:Ek}.

\begin{figure}[H]
  \centering
  \includegraphics[width=0.48\linewidth,height=0.24\linewidth]{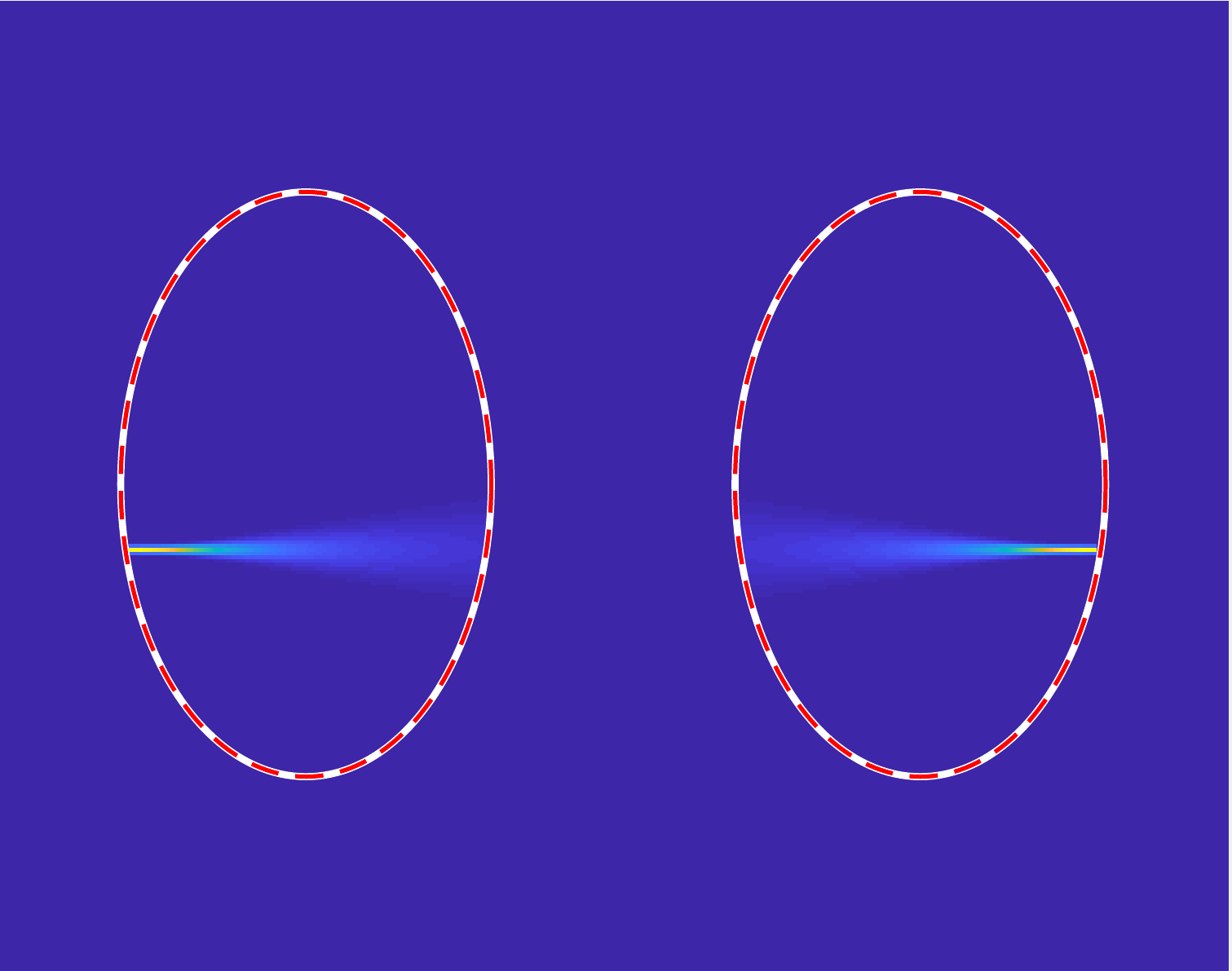}
  \includegraphics[width=0.48\linewidth,height=0.24\linewidth]{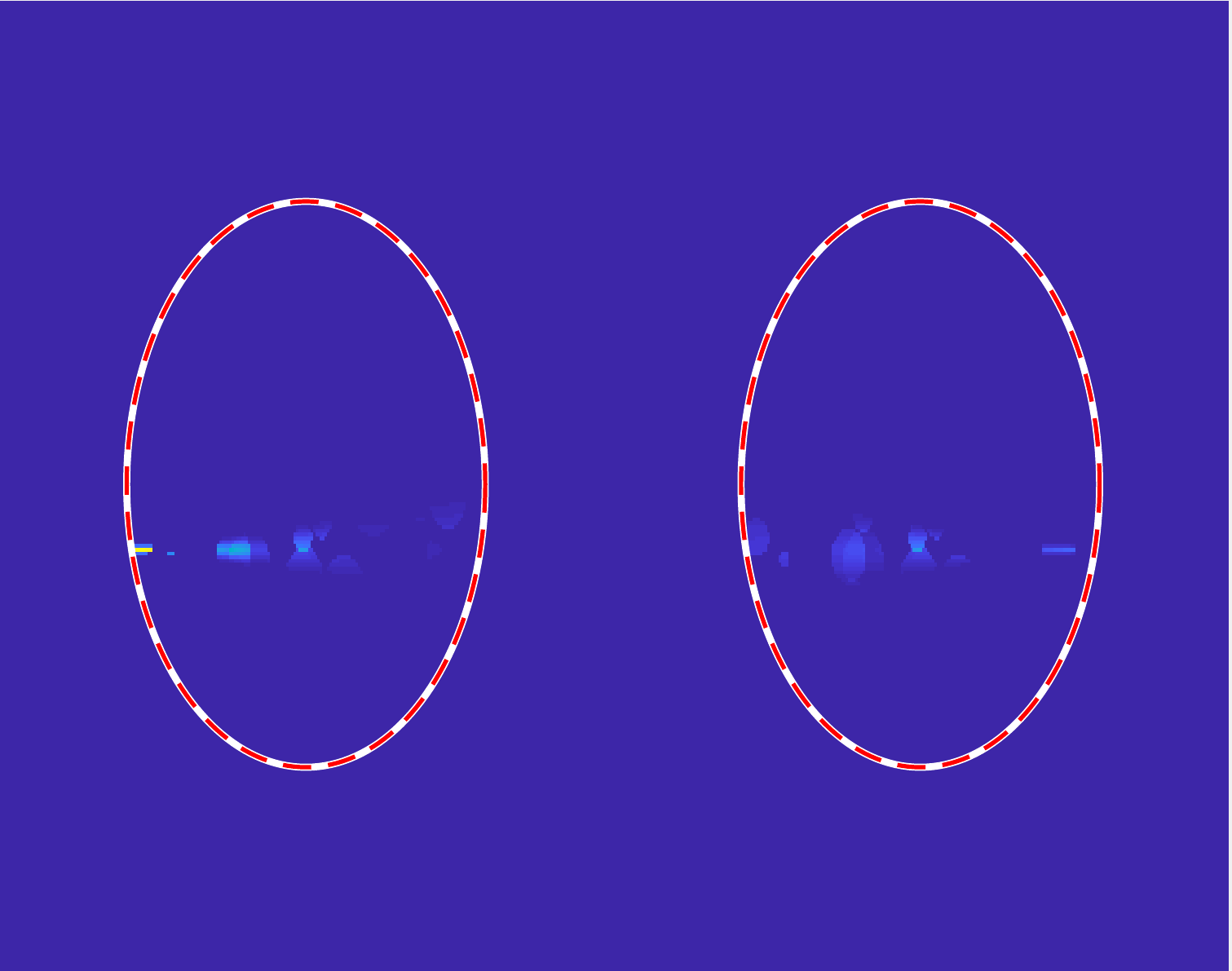}
  \caption{The left image corresponds to the $v_h$ image after illuminating at $h=-0.1406$ from left and right, respectively. In the second image, we show the function $w_h$ for the same height. We included the support of our object in broken red lines for visualization purposes.}
  \label{fig:vh}
 \end{figure}

 Figure~\ref{fig:vh} presents a single simulation of $v_h(x,y)$ when the excitation point occurs at $h=-0.1406$, from both directions (left and right). We also included a visualization of the function $w_h$.

 To achieve the discretization of the equation~\eqref{eq:measure_h}, we define the set of discrete values of $h$ as $\{h_l\}$ for $l=1,\ldots, 2N$ and analogously, for $s$ we consider $\{s_k\}$ for $k=1,\ldots,N.$
 Additionally, as images $a$ and $\mu$ are seen as matrices, we index them as $a_{ij}$ and $\mu_{ij}$ for $i,j = 1,\ldots, N$. Finally, a line of observation is defined by the distance $s_k$, and we denote it by $L_k$.

\begin{figure}[ht]
 \centering
 \includegraphics[scale=1.3]{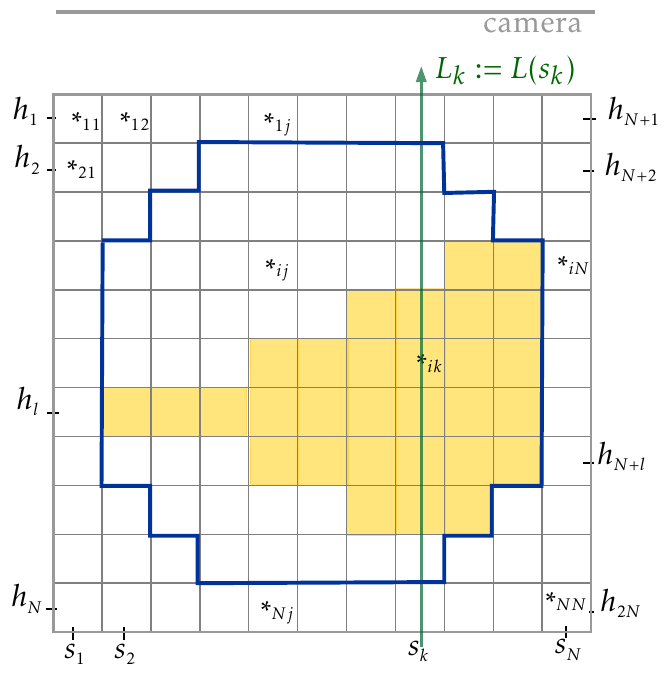}
 \caption{Discretization of the image and the variables used in the AtRt.}\label{fig:new_discretization}
\end{figure}

 In Figure~\ref{fig:new_discretization}, we describe all the discrete variables that we have introduced. The filled pixels represent an example of the discretized function $v_h$ when the excitation occurs at the point $h_l$ of our discrete domain. We denote by $v_{ijl}$ the value of $v_{h}$ in the pixel indexed by $(i,j)$ when $h=h_l$.
 We use the Kronecker delta to determine if a line $L_k$ is intersecting a pixel $(i,j)$, this happens when we are at pixels where $j=k$, then:
 \[
 \delta_{jk} = \left\{
 \begin{array}{ll}
 1,& \text{if $j=k$,}\\
 0,& \text{otherwise.}
 \end{array}
 \right.
 \]
 Then $\mathcal P[\mu] (s_k, h_l)= p_{h_l}(s_k)$ is calculated as:
\begin{eqnarray}
\label{eq:02} \mathcal P[\mu](s_k, h_{l})
& = &  c\sum_{i,j=1}^{N} \delta_{jk}\mu_{ij}v_{ijl}\exp\left(-D_{ik}(a)\right),\\
\label{eq:03} & = &  c\sum_{i=1}^{N} \mu_{ik}v_{ikl}\exp\left(-D_{ik}(a)\right),
\end{eqnarray}
where
\[
D_{ik}(a) = \sum_{z=1}^{i} a_{zk}
\]
is interpreted as partial sums along the columns of the attenuation $a$.

Under this discretization, our set of measurements is of size $2N^2$, for all $(s_k, h_l)$ with two--side excitations (we highlight that the density $\mu$ has $N^2$ pixels that is the amount of unknowns of our problem). In Figure~\ref{fig:measurements}, the first two images represent the matrix of measurements obtained from left and right excitations, respectively. In the third one, the fused image (as in \cite{huisken2012slicing} is presented to compare it with the reconstruction obtained by the proposed model.

\begin{figure}[ht]
  \centering
  \includegraphics[width=0.7\linewidth]{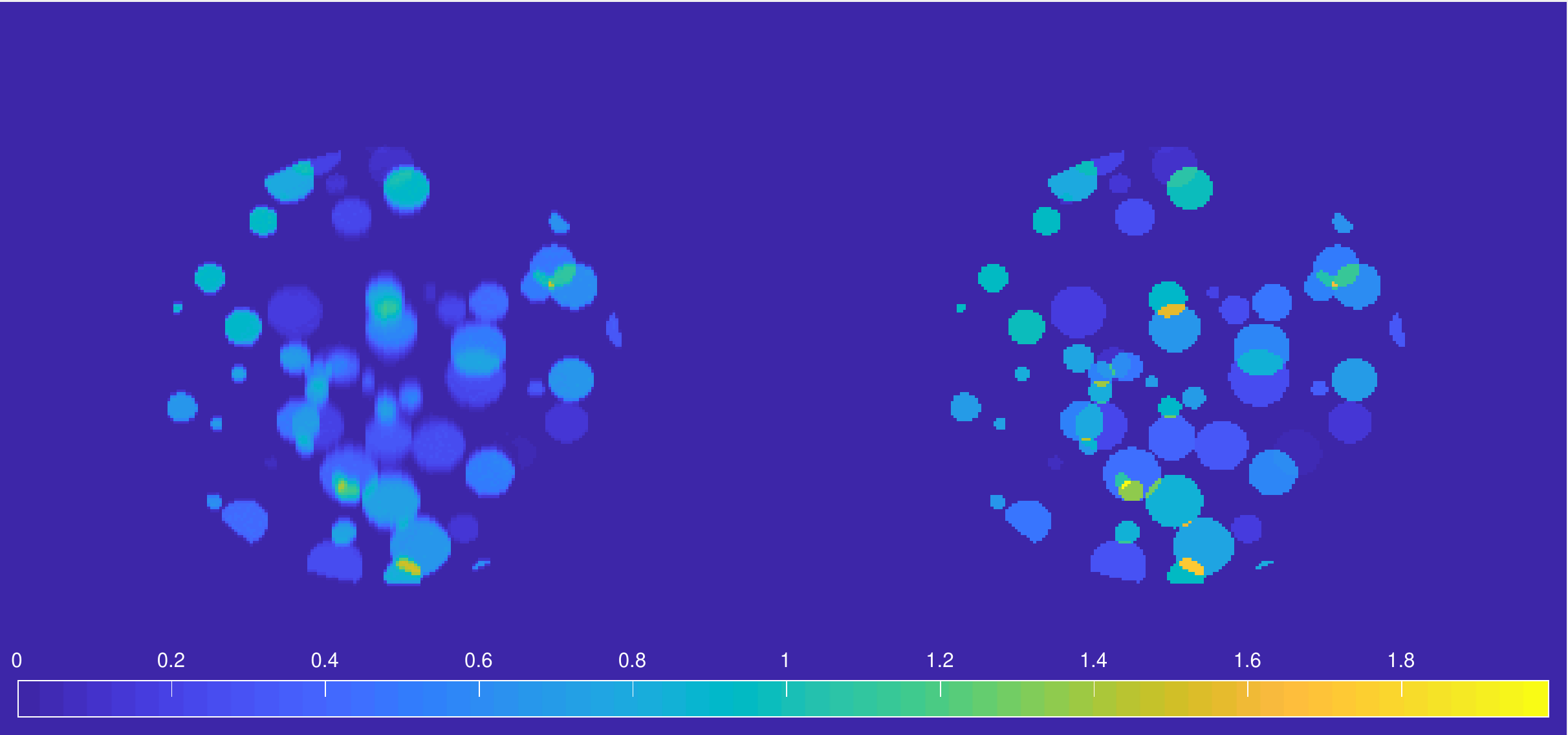}
  \caption{The fused image (see Figure~\ref{fig:measurements}) of measurements (left) compared to the ground truth density (right).}\label{fig:measurementsVSsource}
 \end{figure}

 In Figure~\ref{fig:measurementsVSsource}, we compare the fused image and the ground truth density $\mu$ under the same scale of values. This figure shows that the density that is measured by the camera is not as good and need to be corrected in the central zone, which was our initial motivation. In the next section, we study the numerical inversion of the proposed inverse problem and
 present possible improvements that can be obtained through our approach.

 \subsection{Inverse model}

 We take advantage of the linearity of the operator $\mathcal P$ described in Definition~\ref{def:measurement_operator}, to represent the solution of our discretized inverse problem as the solution of a linear system of the form:
 \begin{equation}\label{eq:linear_system}
 A{\bm\mu}=b,\qquad A\in\Rbb^{m\times n}, \ b\in \Rbb^m,\ \bm\mu\in\Rbb^n.
 \end{equation}

 To build the matrix $A$ associated to our problem, we have to do small changes to the previous discretization. We just reorder $(\mu_{ij})$ as a vector $\bm \mu$ of size $N^2\times 1$ as shown in the expression below. We use the variable $z$ to index pixels, so $z=1, \ldots, N^2$. The same is needed for ${\bm v}_l := (v_{ijl})$:
 \[
 \bm \mu = (\mu_z) = \begin{bmatrix}
 \mu_{11}\\ \mu_{21}\\ \mu_{31}\\ \vdots\\ \mu_{NN}
 \end{bmatrix},\qquad
 {\bm v}_l = (v_{zl})= \begin{bmatrix}
 v_{11l}\\ v_{21l}\\ v_{31l}\\ \vdots\\ v_{NNl}
 \end{bmatrix}, \ \forall l=1,\ldots, 2N.
 \]
 Equivalent to the Kronecker delta we introduce a matrix that can tell us the whole information about the intersections between lines $L_k$ and a pixel $z$. For a fixed pixel $z$ and distance $s_k$, we define
 \[
 w_{zk}= \left\{
 \begin{array}{ll}
 1,& \text{if line $L_k$ crosses the pixel $z$,}\\
 0,& \text{otherwise.}
 \end{array}
 \right.
 \]
 Then, we can write a vector ${\bm w}_k$ of size $N^2\times 1$, as follows:
 \[
 {\bm w}_k = \begin{bmatrix}
 w_{1k},\ w_{2k},\ w_{3k},\  \cdots,\ w_{N^2k}
 \end{bmatrix}^\top.
 \]
 And defining ${\bm W}_{kl} = {\bm v}_l\odot {\bm w}_k$, where $\odot$ represents the Hadamard or point--wise product.
 The only part that needs to be written as a vector in expression~\eqref{eq:02} is the exponential term, for this, we define a matrix ${(D_z)}$ as the cumulative sums of the attenuation matrix $a$ in the direction of the camera. The farther a pixel is from the camera, the greater its accumulated value. As before, we rewrite this matrix as a ($N^2\times 1$)--vector, that we denote by $\bm D$:
 \[
 {\bm D} = \begin{bmatrix}
 D_{1},\ D_{2},\ D_{3},\  \cdots,\ D_{N^2}
 \end{bmatrix}^\top.
 \]
 Now, for each $k$ and $l$, we write a row of our final matrix $A$ as:
 \[
 {\bm a}_{kl} = {\bm W}_{kl}\odot \exp({-\bm D}),
 \]
 where $\exp(-\bm D)$ is understood as the exponential of each component of $\bm D$. Then varying $k$ and $l$, we built $A$ of size $m\times n$, with $m = 2N^2$ and $n=N^2$. To build the vector of measurements $b$, as we obtain our set of observations (as the first two images presented in Figure~\ref{fig:measurements}), we just need to reshape them as a column vector taking row by row and transposing them. The shape of the matrix $A$ and vector $b$ are:
 \begin{align*}
 A &{} = \left[\begin{array}{cccc|cccc|c|cccc}
 {\bm a}_{11} & {\bm a}_{21}& \cdots&{\bm a}_{N1}&
 {\bm a}_{12}& {\bm a}_{22}& \cdots& {\bm a}_{N2}&
 \cdots&
 {\bm a}_{1,2N}& {\bm a}_{2,2N}& \cdots& {\bm a}_{N,2N}
 \end{array}\right]^\top, \\[8pt]
 b &{} = \left[\begin{array}{cccc|cccc|c|cccc}
 {b}_{11}& {b}_{21}& \cdots& {b}_{N1}&
 {b}_{12}& {b}_{22}& \cdots& {b}_{N2}&
 \cdots&
 {b}_{1,2N}& {b}_{2,2N}& \cdots& {b}_{N,2N}
 \end{array}\right]^\top.
\end{align*}

\subsubsection{Solution of the linear system.}

 As the matrix $A$ is sparse and large, a factorization process to solve~\eqref{eq:linear_system} could be impossible or computationally expensive. For this reason, the use of iterative methods is highly desirable to solve this type of linear systems.

 Additionally, we consider that our measurements (represented by the right-hand vector $b$) are corrupted by unknown vector of noise $\varepsilon\in \Rbb^m$, as is usual in the real cases. For the different iterative algorithms that we will present, we assume that at least the norm $\delta := \|\varepsilon\|$ is known.

 Then, due to the ill--posedness produced by the presence of noise and the possible ill-conditioned matrix $A$, a \emph{regularization process} can be used to overcome these issues~\cite{calvetti2004non}.

 The regularized minimization problem associated to the solution of the linear system~\eqref{eq:linear_system} is:
 \begin{equation}\label{eq:regularization_problem}
 \mu = \argmin_{x\in \Rbb^m} \left\{\frac 12 \|Ax-b\|_2^2 + \lambda \mathcal R(x)\right\}
 \end{equation}
 where the data--fit term $\|Ax-b\|_2^2$ forces the problem to find $x$ that remains close to the given data $b$, and the regularizer term $\mathcal R$ is chosen to overcome the particular requirements of each problem. An alternative way to include the regularization is to apply an iterative method directly on the data--fit term and use the number of iterations as stop criteria when semi-convergence is achieved. The general principle of the semi--convergence is to obtain a desired approximation before the noise starts to show up in the current solution~\cite[Chapter 6]{hansen2010discrete}. The algorithms used to solve our problem consider these two possible approaches.

 In the next section, we briefly describe the algorithms that are used to solve our linear system and hence, the inverse problem. We have implemented the discretization of our problem in \textsc{Matlab} and we solve the linear system using the \textsc{IR tools} which are detailed in~\cite{gazzola2019ir}.

 \section{Numerical results}\label{sec:numerical_results}

 In this part, we propose to solve our discrete inverse problem using two different minimization approaches, that we denote by (P1) and (P2) and are defining as follows:
 \begin{align}
 &\tag{P1}\label{eq:P1} \left\{\begin{aligned}
 & \underset{x}{\text{minimize}}
 & & \|Ax -b\|^2_2 \\
 & \text{subject to}
 & & x \in \mathcal C
 \end{aligned}\right.\\[10pt]
 &\tag{P2}\label{eq:P2} \left\{\begin{aligned}
 & \underset{x}{\text{minimize}}
 & & \|Ax -b\|^2_2 +\lambda \TV(x) \\
 & \text{subject to}
 & & x \geq 0
 \end{aligned}\right.
 \end{align}
 The Problem~\eqref{eq:P1} is related to the \emph{semi--convergence} case, where the regularization will be included within the iterations of the optimization algorithms. We will compare the results obtained by five different algorithms: the \emph{Modified residual norm steepest descent method}~\cite{nagy2000enforcing} (\texttt{mrnsd}), the \emph{Flexible CGLS method}~\cite{gazzola2017fast} (\texttt{nnfcgls}), \emph{Simultaneous algebraic reconstruction technique}~\cite{hansen2018air} (\texttt{sart}) and the \emph{Fast Iterative Shrinkage-Thresholding Algorithm} (\texttt{fista})~\cite{beck2009fast} (that solves the Tikhonov problem with box constraints when the parameter $\lambda = 0$, a penalized version is also available if $\lambda \not=0$ but we are not considering this case).


 The problem~\eqref{eq:P2} has the shape of~\eqref{eq:regularization_problem} where we have considered the \emph{total variation} (TV, \cite{rudin1992nonlinear}) as our regularizer $\mathcal R$ .
 To solve it, we use a particular case of the \emph{Projected-restarted iteration method} (PRI)~\cite{calvetti2004non} which incorporates a heuristic TV penalization term~\cite{gazzola2014generalized}. As in~\cite{gazzola2019ir}, we denote this method by (\texttt{htv}).

\subsection{Simulated noise measurements}

To avoid \emph{inverse crime} in our reconstructions, we add noise to our simulated measurements. For this, we consider an scaling factor $\beta$ to generate a poisson distributed noise (since this random variable returns normal values, it is necessary to amplify the signal). The factor $\beta$ controls the level of noise, \emph{i.e.}, if $\beta$ takes large values, we will get lower intensity images and therefore higher poisson noise~\cite{li2017pure}. Accordingly, each pixel value $p$ is replaced by a draw $\beta\cdot \text{Pois}\left(\frac{p}{\beta}\right)$ as in~\cite[eq. 2]{li2017pure}. 

Examples 1 and 2 described below are implemented with values $\beta = 0.01$ and $\beta = 0.001$, respectively.

\subsection{Stop criteria}
In this \textsc{IR tools} package, all algorithms mentioned above used the \emph{discrepancy principle} to stop in the \emph{best} iteration. For the algorithms \texttt{sart}, \texttt{fista}, \texttt{mrnsd} and \texttt{nnfcgls}, this means that the algorithms stop as soon as the relative norm of the residual $b-Ax^{(k)}$ is sufficiently small, typically of the same size as the norm of the noise $\varepsilon$, \emph{i.e.} when
\[
\frac{\|b-Ax^{(k)}\|_2}{\|b\|_2}\leq \eta \cdot \texttt{NoiseLevel}
\]
where $\eta$ is a ``safety factor'' slightly larger than 1, and \texttt{NoiseLevel} is the relative noise $\|\varepsilon\|_2/\|b\|_2$.

For the algorithm \texttt{htv} that is a PRI method with
inner--outer iterations, the discrepancy principle is used to stop the inner iterations, whilst the outer iterations are stopped when $\|x^{(k)}\|$, $\|\TV(x^{(k)})\|_2$ or the value of the regularizer parameter $\lambda$, becomes stable.

\subsection{Initialization}
We use the fused image of measurements (see Figure~\ref{fig:measurements}) as initial value $x^{(0)}$ (see Figure~\ref{fig:measurementsVSsource}), this initializing helps to improve the speed of the algorithms and reduce the number of iterations.

When the parameter $\eta$ is needed, we considered $\eta = 1.01$. Additionally, since we simulate the data as shown in Section~\ref{sec:numerics}, we have at our disposal the true value of the unknown image
$\mathbf{\mu}$ which is included in the algorithm to calculate the relative error. 

\subsection*{Example 1:}
In this first simulated example, we consider that the attenuations $\lambda$ and $a$ are constant over the domain $\Omega$. This means that we are only considering the effects of
the medium where our object of interest in submerged. In Table~\ref{tab:constant_algorithms}, we present the results in terms of number of (outer) iterations, time of execution, the relative error (NRE) and the structural similarity coefficient (SSIM, \cite{wang2004image}) between the reference (true) density and the reconstruction. In this example, all the algorithms present a quantitative improvement compared to the values of the fused image. The \texttt{htv} method gives the smallest NRE value (0.139\%) and \texttt{fista} the highest value of the SSIM (0.98439). In Figure~\ref{fig:constant_zoom}, we can visually compare the different results.

\begin{table}[ht]
\caption{Number of iterations, execution time, relative error and SSIM for the different algorithms when attenuation is assume to be known and constant. The ``fused image'' row corresponds to the third image in Figure~\ref{fig:measurements}, which has been perturbed by noise.}\label{tab:constant_algorithms}
\begin{center}
\begin{tabular}{lcccl}
\bottomrule
Algorithm & iterations & time (s) & $\|x - x^{(k)}\|_2/\|x\|_2$ & SSIM\\\mr
\texttt{fused image} & -- --  &  -- --  &  0.1637   & 0.96402\\
\texttt{fista}       & \,31   & 4.8129   & 0.15077   & $\mathbf{0.98439}$\\
\texttt{htv}         & \,34   & 1.2496   & $\mathbf{0.13914}$   & 0.98349\\
\texttt{mrnsd}       & 150    & 2.8388   & 0.14965   & 0.98278\\
\texttt{nnfcgls}     & 106    & 3.7828   & 0.14001   & 0.98383\\
\texttt{sart}        & \,10   & 1.9969   & 0.15856   & 0.98305\\
\bottomrule
\end{tabular}\\
{\footnotesize $^\ast x$ is the truth solution.}
\end{center}
\end{table}

\begin{figure}[ht]
\begin{center}
\begin{tikzpicture}[spy using outlines={rectangle,white,magnification=2,size=1.5cm, height=5.5cm, connect spies}]
 \node {\pgfimage[height=5cm, width=5cm]{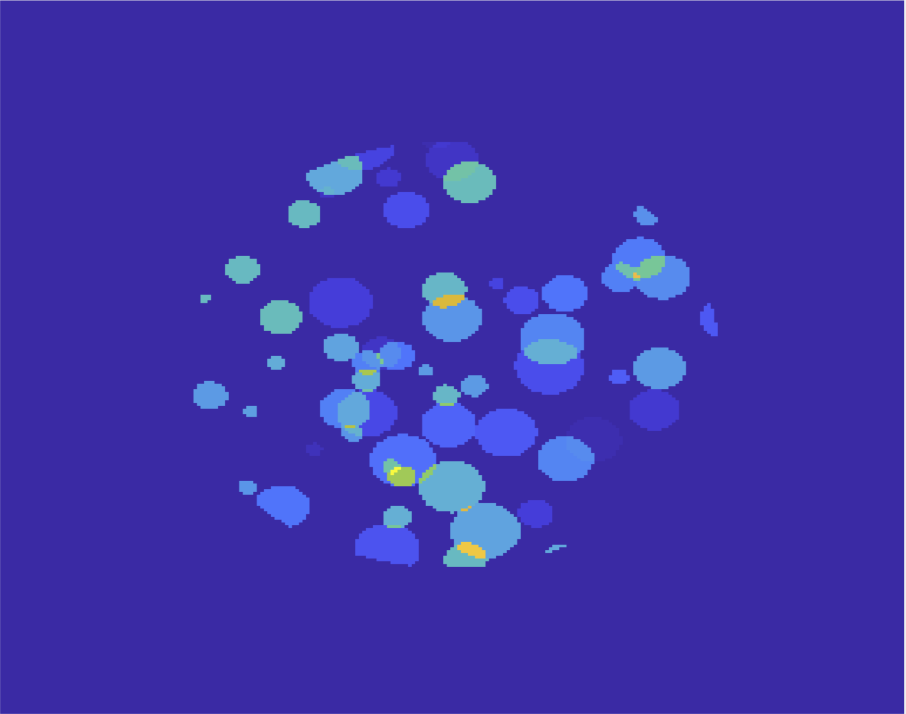}};
 \spy on (0,0) in node [left] at (3.5,0);
 \end{tikzpicture}
\end{center}
 \caption{True simulated density, with zoomed zone to visual comparisons.}\label{fig:zoom_source}
\end{figure}

\begin{figure}[ht]
	\begin{center}
\begin{tabular}{cc}
 \subfloat[Noise measurements]{\begin{tikzpicture}[spy using outlines={rectangle,white,magnification=2,size=1.5cm, height=5.5cm, connect spies}]
 \node {\pgfimage[height=5cm, width=5cm]{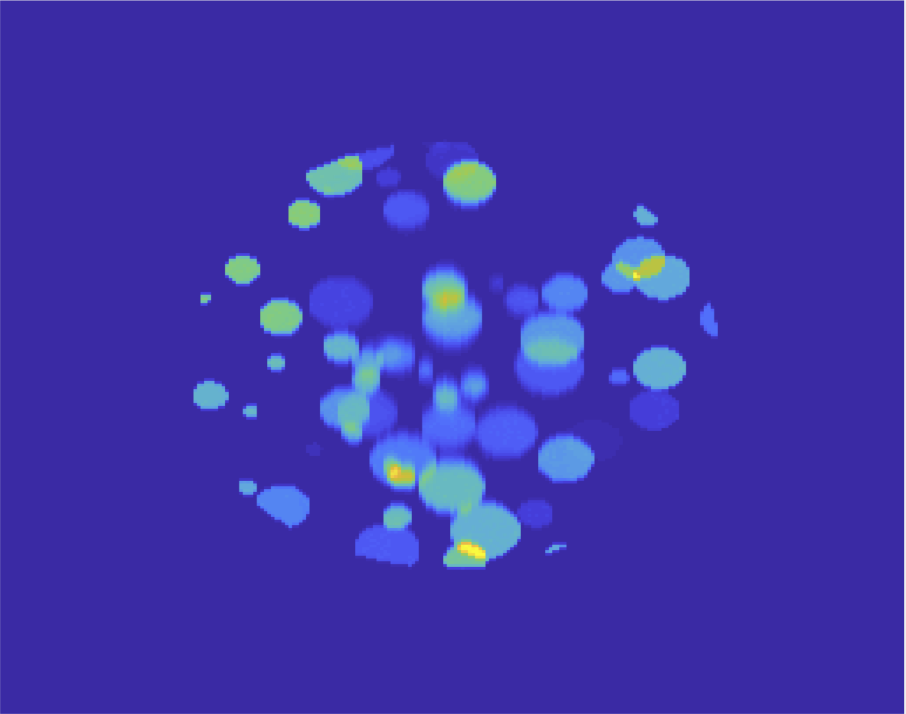}};
 \spy on (0,0) in node [left] at (3.5,0);
 \end{tikzpicture}}&
 \subfloat[fista]{\begin{tikzpicture}[spy using outlines={rectangle,white,magnification=2,size=1.5cm, height=5.5cm, connect spies}]
 \node {\pgfimage[height=5cm, width=5cm]{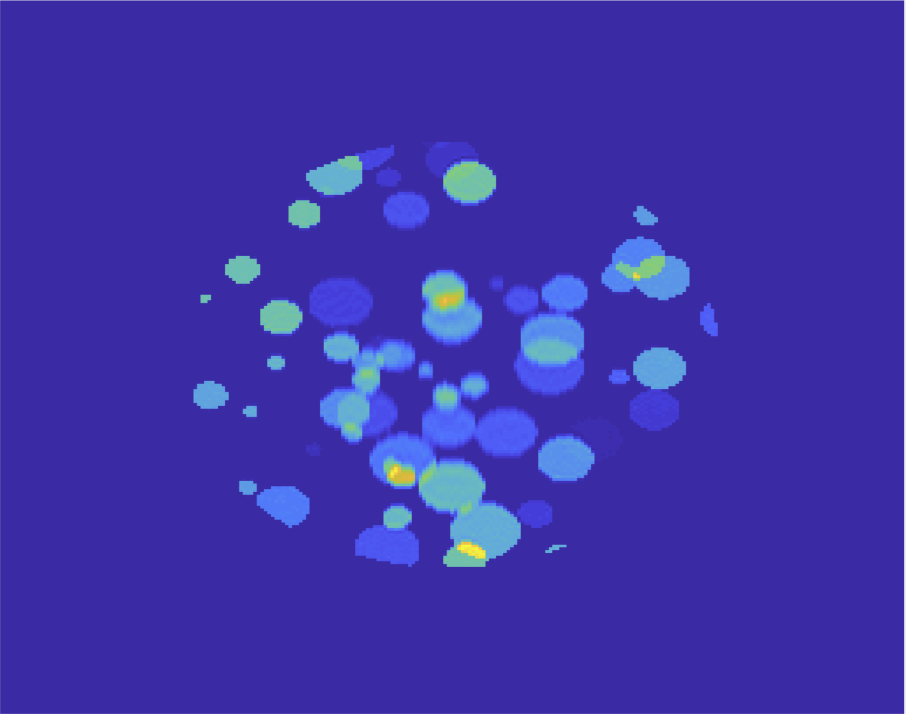}};
\spy on (0,0) in node [left] at (3.5,0);
 \end{tikzpicture}}\\
 \subfloat[htv]{\begin{tikzpicture}[spy using outlines={rectangle,white,magnification=2,size=1.5cm, height=5.5cm, connect spies}]
 \node {\pgfimage[height=5cm, width=5cm]{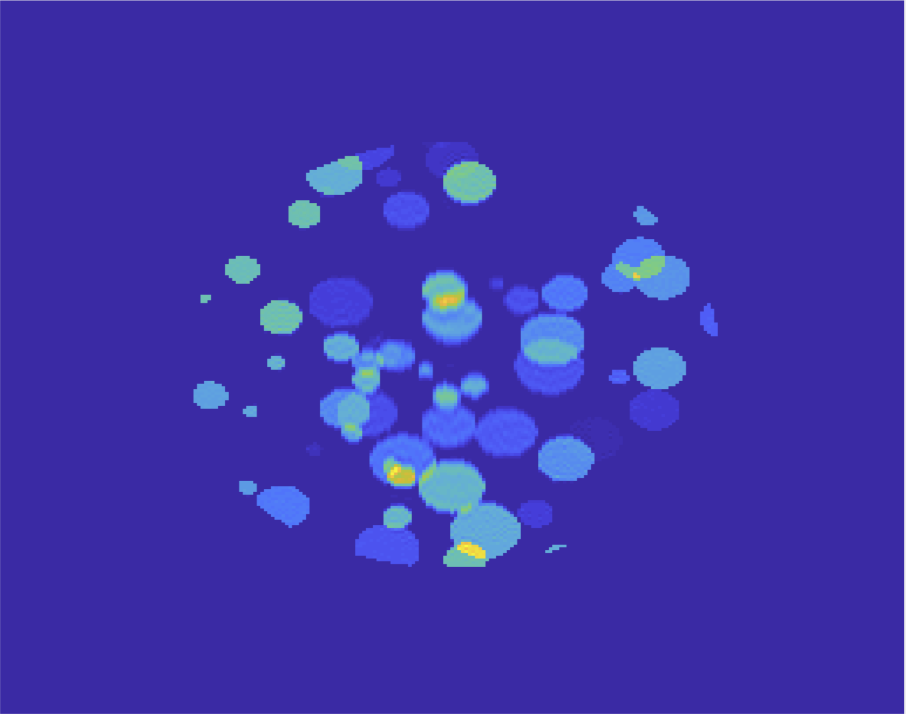}};
\spy on (0,0) in node [left] at (3.5,0);
 \end{tikzpicture}}&
 \subfloat[mrnsd]{\begin{tikzpicture}[spy using outlines={rectangle,white,magnification=2,size=1.5cm, height=5.5cm, connect spies}]
 \node {\pgfimage[height=5cm, width=5cm]{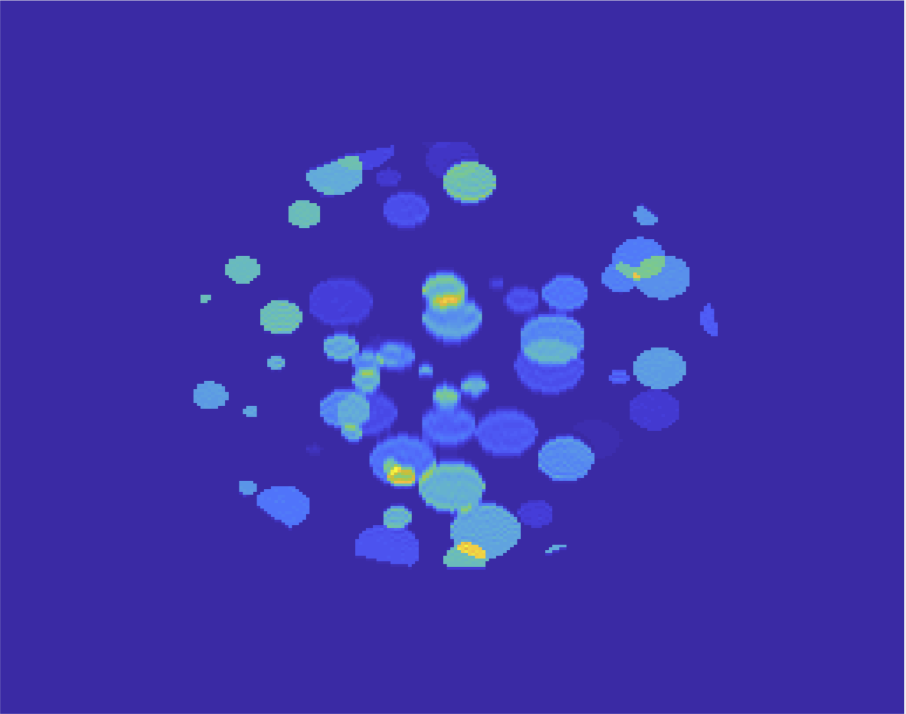}};
\spy on (0,0) in node [left] at (3.5,0);
 \end{tikzpicture}}\\
 \subfloat[nnfcgls]{\begin{tikzpicture}[spy using outlines={rectangle,white,magnification=2,size=1.5cm, height=5.5cm, connect spies}]
 \node {\pgfimage[height=5cm, width=5cm]{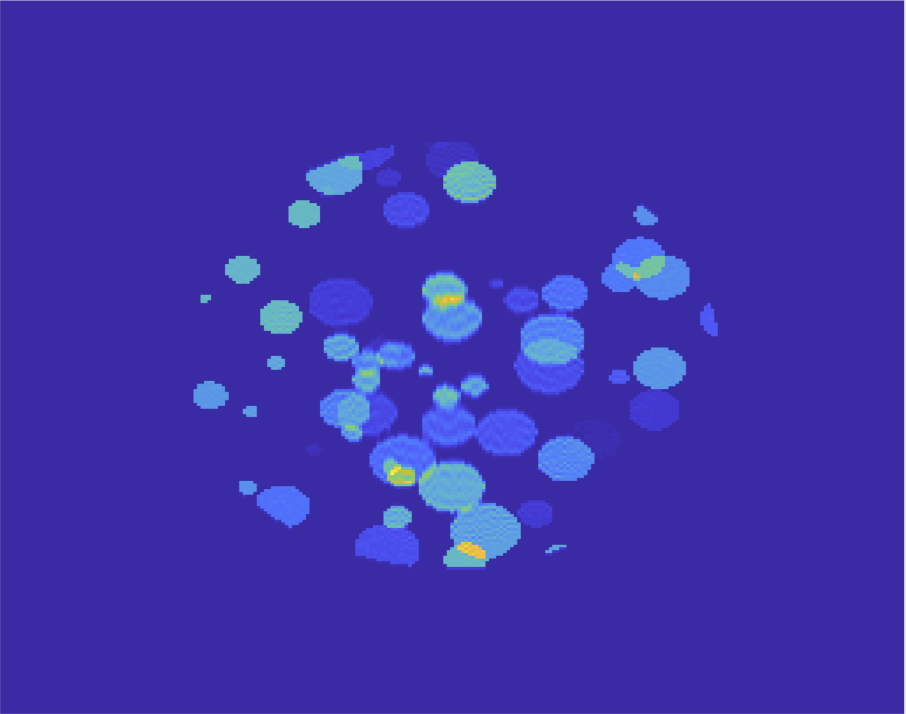}};
\spy on (0,0) in node [left] at (3.5,0);
 \end{tikzpicture}}&
 \subfloat[sart]{\begin{tikzpicture}[spy using outlines={rectangle,white,magnification=2,size=1.5cm, height=5.5cm, connect spies}]
 \node {\pgfimage[height=5cm, width=5cm]{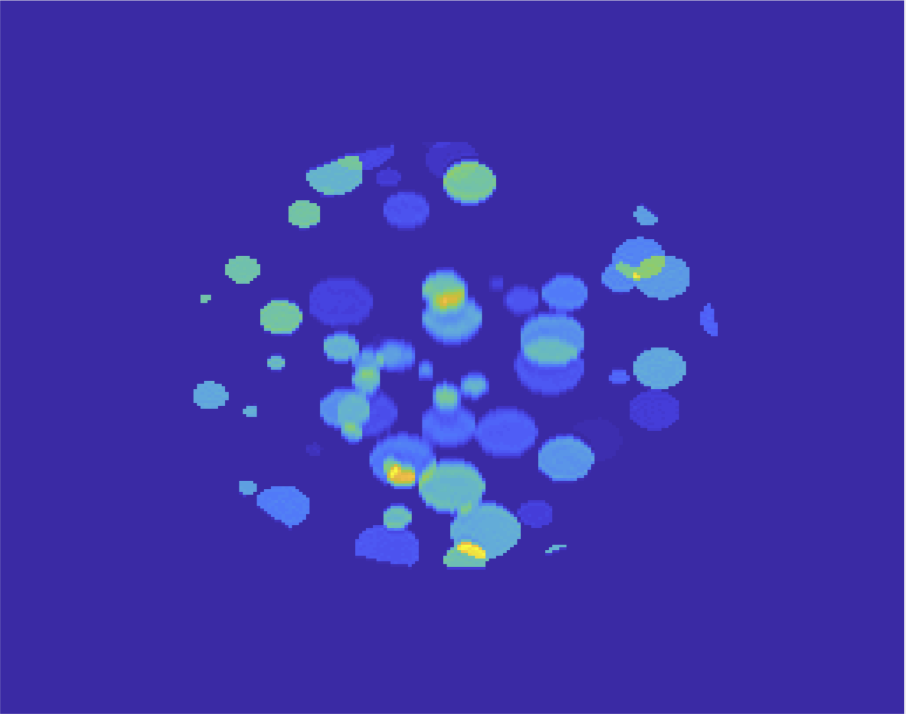}};
\spy on (0,0) in node [left] at (3.5,0);
 \end{tikzpicture}}\\
 \end{tabular}
\end{center}
 \caption{For Example 1: zoomed images to visualize the difference between the reconstructions.}\label{fig:constant_zoom}
 \end{figure}
\begin{figure}[ht]
\begin{center}
 \includegraphics[width=\linewidth]{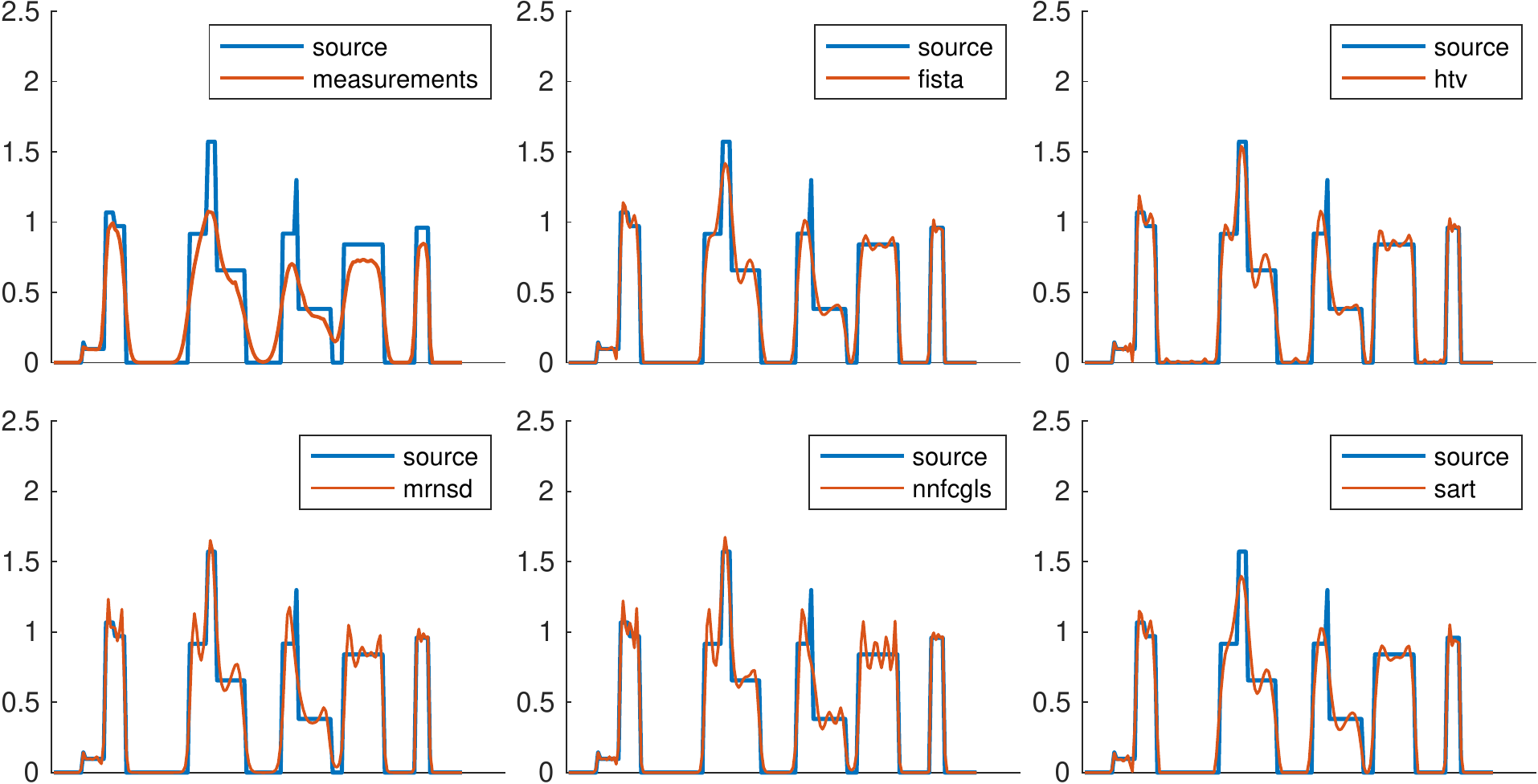}
\end{center}
\caption{For Example 1: Profiles of reconstruction at $x=1$ that corresponds to the column 129 of the images.}
\label{fig:line_comparison_constant}
\end{figure}

In Figure~\ref{fig:line_comparison_constant}, we draw the profiles of the reconstructions along $x=1$ in order to observe the improvements reached in the central region of the image.

\subsection*{Example 2:}
\begin{figure}[ht]
 	\begin{center}
 \begin{tabular}{cc}
  \subfloat[Noise measurements]{\begin{tikzpicture}[spy using outlines={rectangle,white,magnification=2,size=1.5cm, height=5.5cm, connect spies}]
  \node {\pgfimage[height=5cm, width=5cm]{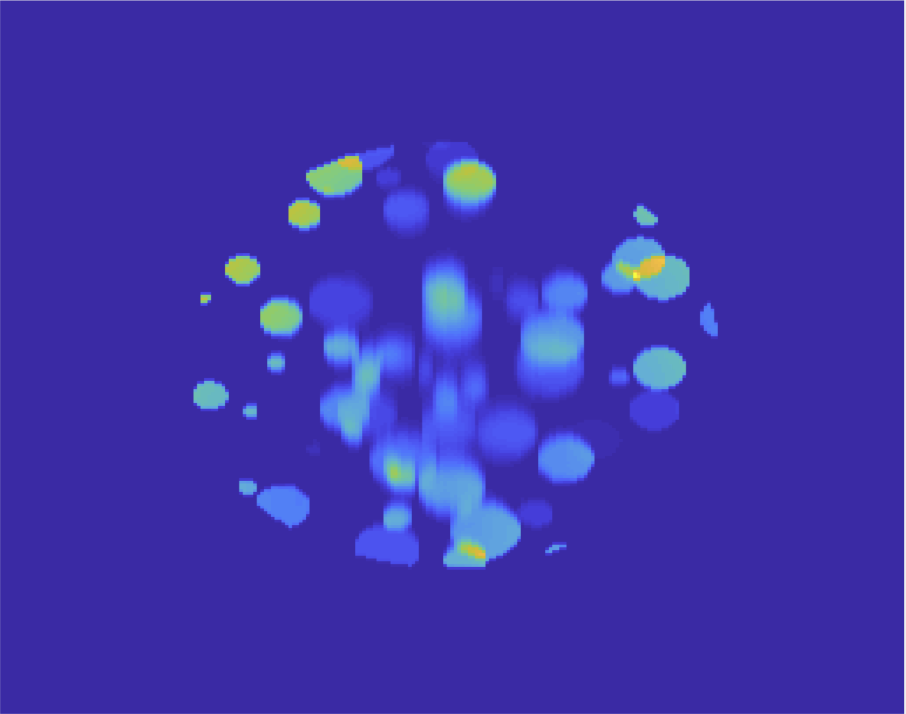}};
  \spy on (0,0) in node [left] at (3.5,0);
  \end{tikzpicture}}&
  \subfloat[fista]{\begin{tikzpicture}[spy using outlines={rectangle,white,magnification=2,size=1.5cm, height=5.5cm, connect spies}]
  \node {\pgfimage[height=5cm, width=5cm]{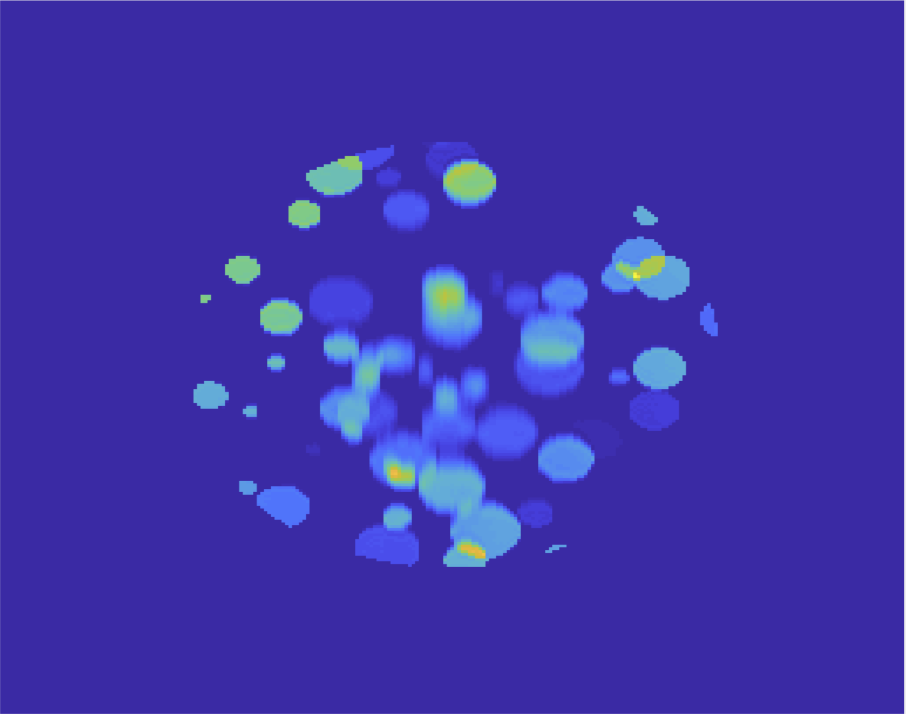}};
 \spy on (0,0) in node [left] at (3.5,0);
  \end{tikzpicture}}\\
  \subfloat[htv]{\begin{tikzpicture}[spy using outlines={rectangle,white,magnification=2,size=1.5cm, height=5.5cm, connect spies}]
  \node {\pgfimage[height=5cm, width=5cm]{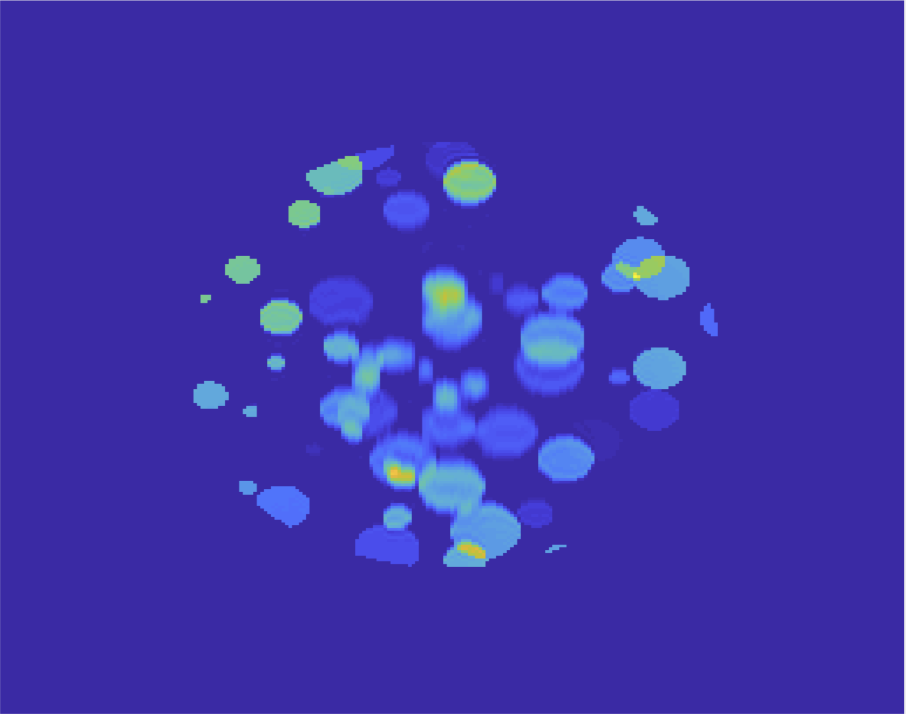}};
 \spy on (0,0) in node [left] at (3.5,0);
  \end{tikzpicture}}&
  \subfloat[mrnsd]{\begin{tikzpicture}[spy using outlines={rectangle,white,magnification=2,size=1.5cm, height=5.5cm, connect spies}]
  \node {\pgfimage[height=5cm, width=5cm]{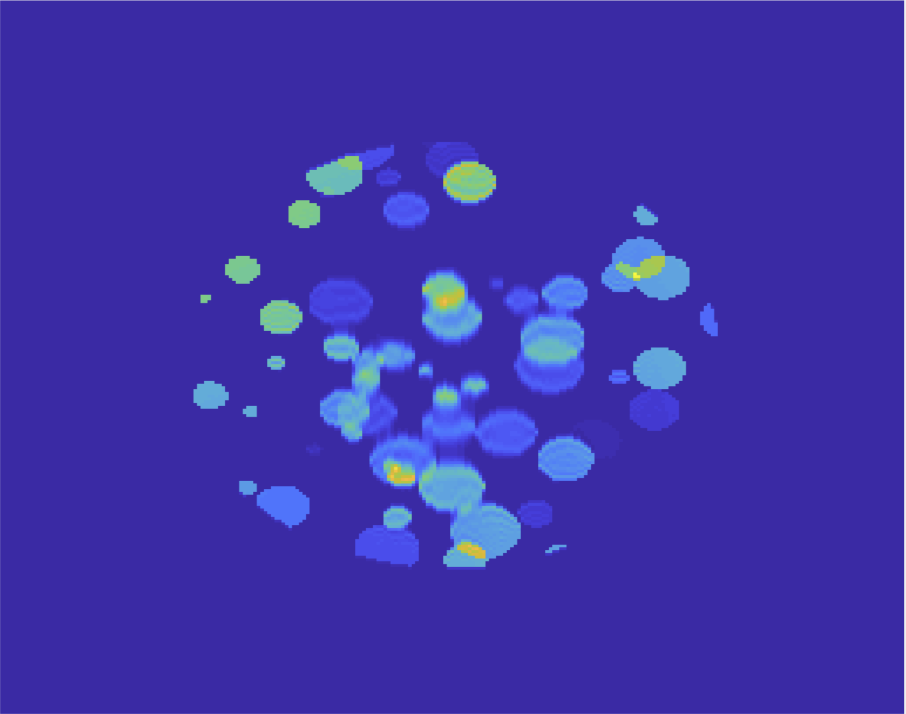}};
 \spy on (0,0) in node [left] at (3.5,0);
  \end{tikzpicture}}\\
  \subfloat[nnfcgls]{\begin{tikzpicture}[spy using outlines={rectangle,white,magnification=2,size=1.5cm, height=5.5cm, connect spies}]
  \node {\pgfimage[height=5cm, width=5cm]{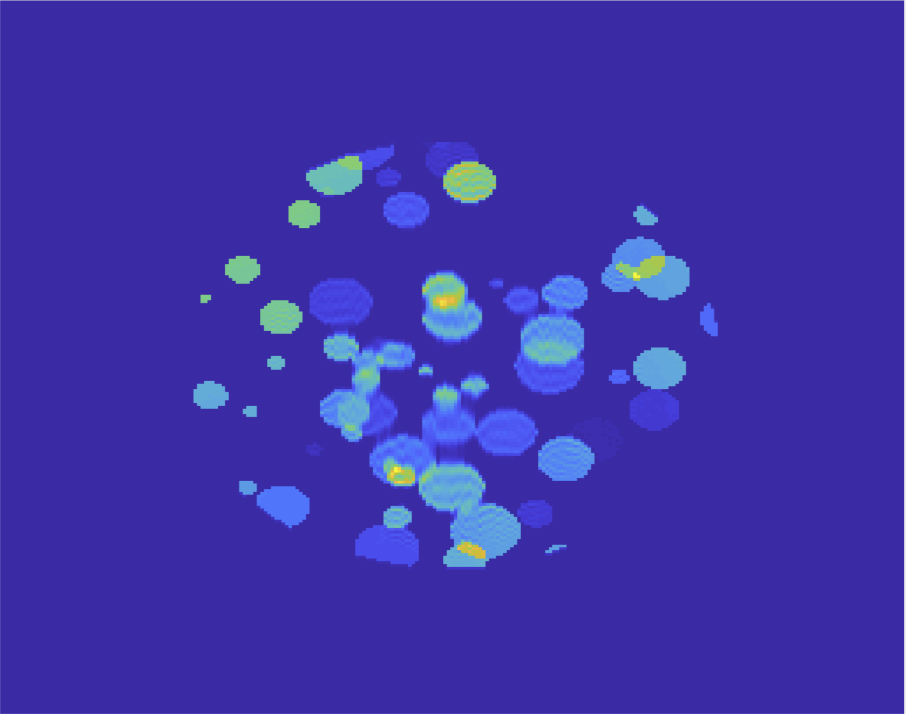}};
 \spy on (0,0) in node [left] at (3.5,0);
  \end{tikzpicture}}&
  \subfloat[sart]{\begin{tikzpicture}[spy using outlines={rectangle,white,magnification=2,size=1.5cm, height=5.5cm, connect spies}]
  \node {\pgfimage[height=5cm, width=5cm]{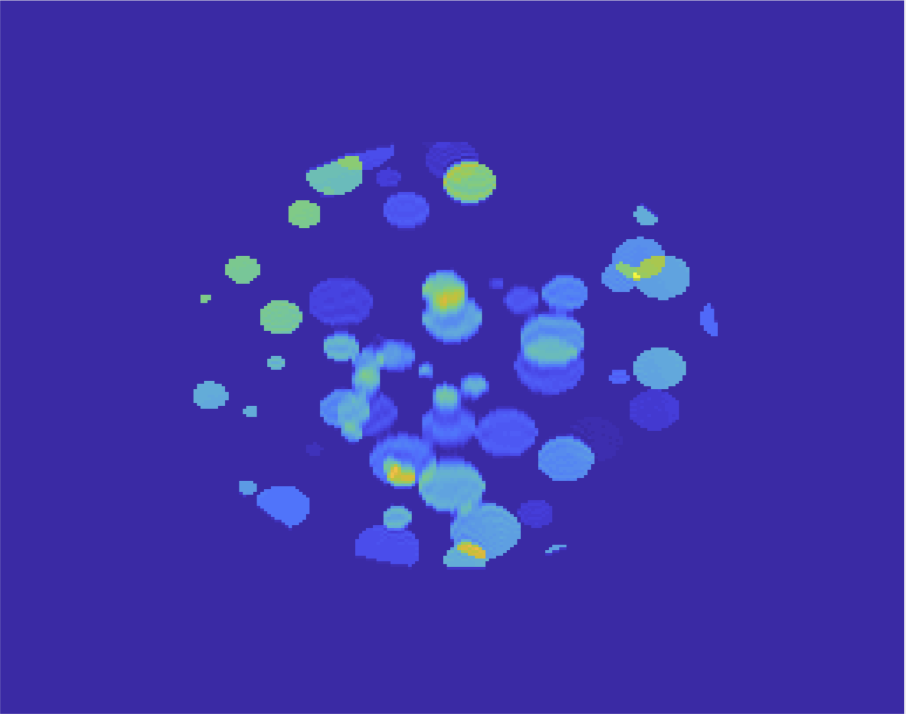}};
 \spy on (0,0) in node [left] at (3.5,0);
  \end{tikzpicture}}\\
  \end{tabular}
 \end{center}
  \caption{For Example 2: zoomed images to visualize the difference between the reconstructions. The images are re-scaled to the range of the ground truth density.}\label{fig:variable_zoom}
  \end{figure}

\begin{figure}[ht]
 \begin{center}
  \includegraphics[width=\linewidth]{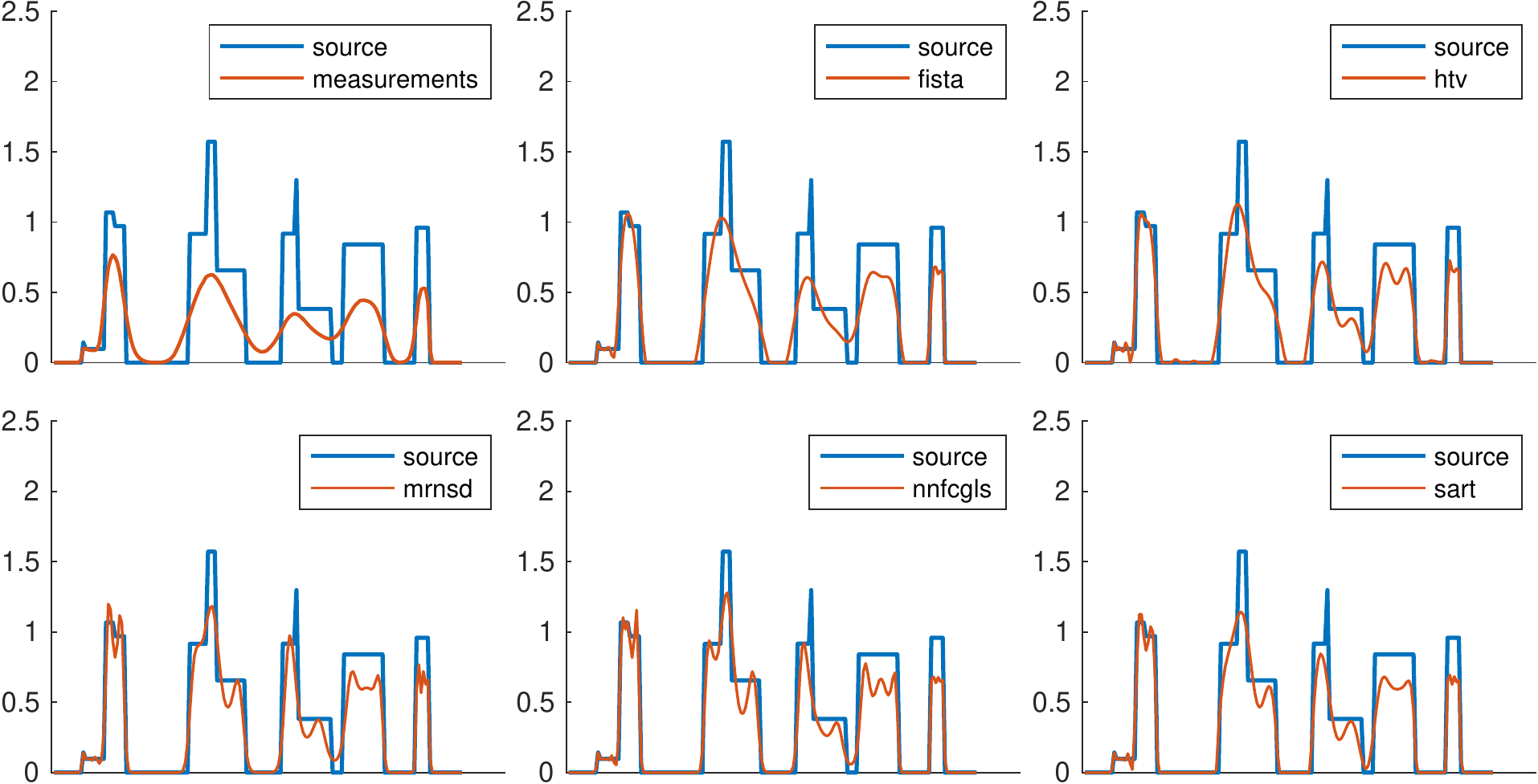}
 \end{center}
 \caption{For Example 2: Profiles of reconstruction at $x=1$ that corresponds to the column 129 of the images.}
 \label{fig:line_comparison_variable}
 \end{figure}

In this case, the simulated measurements are generated using variables attenuations $\lambda$ and $a$, in order to include some attenuation effects produced by the presence of the fluorescent molecules. However, as in more real cases, the attenuation could be also unknown, we propose to reconstruct the density $\mu$ with a constant attenuation $a$ which could be experimentally determined. In our case,
we take $a=1.1$ over $\Omega$. We have included Poisson Noise with $\texttt{NoiseLevel}=0.01$. The results are presented as before in Figures~\ref{fig:variable_zoom}--\ref{fig:line_comparison_variable} and Table~\ref{tab:variable_algorithms}. The values of the \texttt{nnfcgls} and \texttt{sart} methods are slightly better than the other algorithms, but all of them improve the fused image values.

We do not focus on which algorithm is better; we are just interested in the improvements observed in the proposed reconstruction independently of the selection of the optimization algorithm.

\begin{table}[ht]
    \caption{Number of iterations, execution time, relative error and SSIM for the different algorithms when the attenuation is variable but is considered as constant during the reconstruction. The ``fused image'' row corresponds to the third image in Figure~\ref{fig:measurements}, which has been perturbed by poisson noise.}\label{tab:variable_algorithms}
\begin{center}
\begin{tabular}{lcccl}
\bottomrule
 Algorithm & iterations & time (s) & $\|x - x^{(k)}\|_2/\|x\|_2$ & SSIM\\\mr
\texttt{fused image} & -- --  &  -- --   & 0.40466   & 0.92454\\
\texttt{fista}       & \,\,\,29   & 5.3721   & 0.29567   & 0.95875\\
\texttt{htv}         & \,\,\,41   & 2.0133   & 0.26267   & 0.96326\\
\texttt{mrnsd}       &$>2000$ & 27.976   & 0.24255   & 0.97345\\
\texttt{nnfcgls}     &$>2000$ & 96.577   & $\mathbf{0.22798}$   & 0.97695\\
\texttt{sart}        &$>2000$ & 45.802   & 0.22783   & $\mathbf{0.97994}$\\
\bottomrule
\end{tabular}\\[5pt]
{\footnotesize $^\ast x$ is the truth solution, the symbol $>$ means stops with a maximum number of iterations.}
\end{center}
\end{table}

 In Figure~\ref{fig:line_comparison_variable}, we draw the profiles of the reconstructions along $x=1$ as before. Here we observe that the assumption of the attenuation is constant implies in some parts a underestimation of the true value. This will depend directly from the constant value that we choose for $a$.

 \section{Conclusions and outlook}

We presented a novel mathematical model for the Light Sheet Fluorescence Microscopy. To our best knowledge, this is the first approach in this direction and is an initial step in trying to understand and tackle some of the issues observed in LSFM. This work shows that by considering the acquisition of the density $\mu$ as an inverse problem a better reconstruction can be obtained, compared to the traditional merging method that is currently used.

From the theoretical point of view, we presented a uniqueness result for the proposed inverse problem, by reducing it to the recovery of the initial condition in a heat equation with measurements in a space--time curve. The stability in the reconstruction of $\mu$ is not considered in this article. However, due to the clear link between the microscopy inverse problem and backward heat propagation the former is expected to be severely ill-posed. The question then is whether Logarithmic stability is the optimal result or if it is possible to obtain a H\"older-type inequality, this kind of result would also open the door to obtain stability results for more physically complete models. This type of question are expected to be addressed in future works.

Additional future work also includes the extension of these results to the three dimensional case, where some extra assumptions might be necessary and we would need to discuss a light-sheet illumination or a beam illumination as the natural extension of the technique presented here.

Questions about a simultaneous reconstruction are also open. For example, about the possibility of recovering the density and the attenuation (either in the illumination or fluorescence) at the same time, by considering additional measurements when rotating the object in multiple directions.

A more ambitious extension of this work would be to consider more complete and less simplified physics for the illumination and fluorescence stages. In this paper we are heavily reliant in the explicit solution of the Fermi pencil beam equation, which makes it very challenging to extend our results to other illumination models. We are also considering a perfect collimation of the fluorescence measurement and different collimation schemes would give rise to other difficulties. Another ambitious extension of this work would be to include the stochastic nature of the fluorescence stage, which would require an MLEM or similar reconstruction techniques to be considered.

 %
 %
 %

 \section*{Acknowledgments}

E.C. was partially funded by CONICYT-PCHA/Doctorado Nacional/2016-21161721 grant, by SENESCYT/Convocatoria2015 and Project UCH-1566 from the Department of Mathematical Engineering at Universidad de Chile.

\noindent A.O. was partially funded by CONICYT grant Fondecyt \#1191903, CONICYT Basal Program PFB-03 (AFB170001) and MathAmsud 18-MATH-04 and CONICYT/FONDAP/15110009.

\noindent M.C. was partially funded by CONICYT grant Fondecyt \#1191903
and M.C. thanks Bo\u{g}azi\c{c}i University, Istanbul, Turkey, as part of this work was completed as a visiting researcher at the institution.

\noindent S.H. and V.C. are part of SCIAN-Lab funded by Fondecyt \#1181823, EQM140119, CONICYT (PIA ACT 1402), CENS CORFO (16CTTS-66390) and BNI (ICM P09-015-F). SCIAN-Lab is a selected member of the German-Chilean Center of Excellence Initiative (DAAD 57220037 and 57168868). V.C. is also partially funded by CONICYT grant
Fondecyt \#11170475.

\noindent B.P. was partially funded by ONR grant N00014-17-1-2096.

\noindent We acknowledge M.D. Miguel Concha for providing us with light-sheet microscopy data (funded by Fondequip EQM130051).

 \section*{References}
 \bibliographystyle{plain}
 \bibliography{revised_article_arXiv_v2}

 \end{document}